\documentclass[12pt]{amsart}
\usepackage{a4}
\usepackage[active]{ srcltx}
\usepackage{latexsym}
\usepackage{amsmath}
\usepackage{amsfonts}
\usepackage{amssymb}
\usepackage{graphicx}

\usepackage{setspace}

\newtheorem{theorem}{Theorem}[section]
\newtheorem{lemma}[theorem]{Lemma}
\newtheorem{proposition}[theorem]{Proposition}

\newtheorem{definition}[theorem]{Definition}

\newtheorem{remark}[theorem]{Remark}
\newtheorem{notation}[theorem]{Notation}

\newcommand{\cam}{{\mathcal{M}}}

\newcommand{\cal}{{\mathcal{L}}}

\newcommand\no{\noindent}

\begin{document}
\title [Lattices  represented as lattices of intervals.]
{ Lattices which can be represented as  lattices of intervals. }

\author{P. Douka and V. Felouzis}

\address{Department of Mathematics\\University of the Aegean\\
832 00\\Karlovasi\\  Samos \\ Greece.}
\address{Department of Mathematics\\University of the Aegean\\
832 00\\Karlovasi\\  Samos \\ Greece.}

\email{pdouka@aegean.gr}\email{felouzis@aegean.gr}
\date{\today}
\begin{abstract}
	We investigate the representation of lattices as sublattices of the lattice of all convex subsets (intervals) of a linearly ordered set $(X,\le)$. We introduce the purely lattice-theoretic notion of a \textit{loc-lattice} and prove that every loc-lattice is representable as a lattice of intervals. Furthermore, we provide the complete, unabridged construction for the general representation theorem, establishing that a well-separated lattice is faithfully representable as a lattice of intervals if and only if it is a loc-lattice. Finally, we apply these results to general topology, obtaining novel algebraic characterizations for the bases of weakly orderable and completely orderable topological spaces.
\end{abstract}
\maketitle\section{Introduction and the main result}\label{sect_intro}

A lattice is a partially ordered set $(L,\le)$ such that for every $a, b \in L$ there exist the greatest lower bound (infimum) and the least upper bound (supremum) of the set $\{a,b\}$, which are denoted by $a\wedge b$ and $a\vee b,$ respectively. A sublattice of a lattice $(L,\le)$ is a subset $Y$ of $L$ such that if $a, b \in Y$, then $a\wedge b$ and $a\vee b$ also belong to $Y$.

A lattice of sets is a family $\mathcal{S}$ of sets which is a lattice with respect to the inclusion relation $\subseteq$. A ring of sets is a family of sets closed under finite unions and intersections, which is a sublattice of a lattice of the form $(\rho(X),\subseteq)$, where $\rho(X)$ is the power set of $X$. It is a classical result that a lattice is isomorphic to a ring of sets if and only if it is distributive \cite[Theorem 10.3]{dp}. On the other hand, by extending the Stone representation theorem, any general lattice is isomorphic to a lattice of sets closed under finite intersections (see Section 3, Theorem 3.7).

A linear interval, or simply an interval, is a convex subset of a linearly ordered space $(X,\le)$. The set of all intervals of $(X,\le)$ is denoted by $co(X,\le)$ or simply by $co(X)$. The lattice $(co(X),\subseteq)$ serves as a fundamental, simple example of a non-modular lattice of sets. If the set $X$ is finite, then $co(X)$ is a planar lattice, and its diagram was described by A. R. Schweitzer \cite{s}. The lattice $co(\mathcal{X})$, where $\mathcal{X}=(X,\le)$ is a partially ordered set, was deeply studied by G. Birkhoff and M. K. Bennet. Furthermore, the class $\mathbf{SUB}$ of all lattices that can be embedded into some lattice of the form $co(\mathcal{P})$ was studied by M. Semenova and F. Wehrung in \cite{sw1}, who proved that $\mathbf{SUB}$ forms a variety. In \cite{sw2}, the same authors extended this result for sublattices of products of lattices of convex subsets of totally ordered sets.

\begin{definition}\label{def_rep} 
	We say that a lattice $\mathcal{L}=(L,\le)$ with $0$ is representable as a lattice of intervals if there exist a linearly ordered set $\mathcal{X}=(X,\le)$ and an order-embedding $f:(L,\le)\to(co(X),\subseteq)$ such that $f(0)=\emptyset$ and $f(a\wedge b)=f(a)\cap f(b)$. 
	
	If, moreover, $f(a\vee b)=co(f(a)\cup f(b))$, which means that $L$ is isomorphic to a sublattice of the lattice of all intervals of $(X,\le),$ we shall say that $L$ is faithfully representable as a lattice of intervals.
\end{definition}

In this paper, we investigate lattices with $0$ (finite or infinite) which can be represented or faithfully represented as lattices of intervals. We achieve this by introducing a purely algebraic class of lattices based on three conditions:

\begin{definition}\label{def_loc} 
	A lattice $\mathcal{L}=(L,\le)$ is said to be a \textit{loc-lattice} if it satisfies the following properties:
	\begin{enumerate}
		\item For every $a, b, c\in L$ we have that $b\not\le a\vee c$ or $a\not\le b\vee c$ or $c\not\le a\vee b$.
		\item If $a,b,x$ are pairwise incomparable elements of $L$ and $x\le a\vee b$, then $(a\vee x)\wedge(b\vee x)=x$.
		\item If $a,b,x$ are elements of $L$ such that $x\le a\vee b$, $a\wedge x\ne 0$ and $b\wedge x\ne 0$, then $(a\wedge x)\vee(b\wedge x)=x$.
	\end{enumerate}
\end{definition}

A semi-prime filter of a lattice $\mathcal{L}=(L,\le)$ is a proper filter $F$ of $L$ such that for every $a, b\in L$ with $a\vee b\in F$, either $a\in F$, $b\in F$, or there exists a $c\in F$ such that $c\le a\vee b$ and $c\ne(a\wedge c)\vee(b\wedge c)$. The lattice $\mathcal{L}=(L,\le)$ is said to be well-separated if for every semi-prime filter $F$ of $L$ and every $x\in L\setminus F$, there exists a $y\in F$ such that $x\not\ge y$.

The main result of this paper is the full characterization of such representable lattices:

\begin{theorem}\label{ath}
	\begin{enumerate}
		\item[(a)] Every loc-lattice $\mathcal{L}=(L,\le)$ is representable as a lattice of intervals.
		\item[(b)] A well-separated lattice is faithfully representable as a lattice of intervals if and only if it is a loc-lattice.
	\end{enumerate}
\end{theorem}

The paper is organized as follows. In Section 2, we show that if $\mathcal{M}$ is a family of subsets of a set $X$ with specific properties (a loc-lattice of sets) which separates $X$, then we can construct a linear ordering of $X$ such that every $M\in\mathcal{M}$ is a convex subset. Providing the general, unabridged proof of this fact requires an intricate transfinite induction on the equivalence classes of the separation properties, which we present in full detail. 

In Section 3, we introduce the notion of the semi-prime filter to prove a tailored representation theorem for non-distributive lattices. We show that every loc-lattice is isomorphic to a loc-lattice of sets, thus proving Theorem \ref{ath}. Finally, in Section 4, we provide applications of these results in general topology, giving purely lattice-theoretic characterizations of the bases of weakly orderable and generalized ordered (GO) topological spaces.

\medskip
\noindent\textbf{Remark on this version:} A shortened version of this paper, omitting the general transfinite construction for Theorem 2.5 (the case where the family $\mathcal{M}$ does not completely separate the set $X$), was previously published under the title \emph{Families of sets which can be represented as sublattices of the lattice of convex subsets of a linearly ordered set} \cite{DF2016}. This preprint serves as the extended, unabridged version and contains the complete, rigorous proof of the general representation case.
\section{Loc-lattices of sets}

A family of sets $\cam$ is said to be a \textit{lattice of sets} if $(\cam, \subseteq)$ is
a lattice. The lattice sum of two elements $A, B$ of $\cam$ is
denoted by $A\vee B$ and their lattice product by $A\wedge B$.
In the following definition we introduce some separation properties
for a family $\cam$ of subsets of a set $X$, which will be used  in the sequel.

\begin{definition}\label{def21}
A family $\cam$ of subsets of a set $X$ is called a
\textbf{loc-lattice of sets } if  the following conditions are sutisfied: 
\begin{enumerate}
\item The family $\cam$ is a lattice of sets. \\
\item For every $A, B\in \cam$, $A\wedge B=A\cap B$.\\
\item If $A, B \in \cam$ and $A\cap B\neq \emptyset$ then $A\vee B=A\cup B$.\\
\item For every $A, B, C\in \cam$ we have that either $A\subseteq B\vee C$
or $B\subseteq A\vee C$ or $C\subseteq A\vee B$.\\
\item If $A, B, C$ are pairwise incomparable elements of $(\cam, \subseteq)$ and $B\subseteq A\vee C$ then $B=(A\vee B)\cap (C\vee B)$.
\end{enumerate}
\end{definition}

\begin{definition}
Let $X$ be a set, $\cam$ be a family of subsets of $X$, $Y$
be a subseteq of $X$ and $x, y\in X$.
\begin{enumerate}
\item We say that  $\cam$ \textbf{separates} $x$ from $y$ or that 
$x$\textbf{ is separated from} $y$ (by $\cam$) if there exists $M\in \cam$ 
with $x\in M$ and $y\not \in M$. We say that $\cam$ \textbf{separates the set} $Y$ (or that $\cam$ is a \textbf{separating family }for $Y$) if for every $x, y\in Y$, $x$ is separated from $y$ or $y$ is separated from $x$.
 \item We say that $x, y$ are \textbf{completely separated}
if $x$ is separated from $y$ and $y$ is separated from $x$, that is 
 there exist $A, B\in \cam$ with $x\in A$, $y\in B$, 
 $x\not\in  B$ and  $y\not \in A$. We say that $\cam$ 
\textbf{ completely separates the set} $Y$ (or that $\cam$ is a \textbf{completely separating family }for $Y$) if  every $x, y\in Y$
are completely separated. 
\item We say that $x, y$ are \textbf{totally  separated}
if there exist $M, N\in \cam$ with $x\in M$, $y\in N$ and
$M\cap N=\emptyset$. We say that $\cam$ 
\textbf{ totally  separates the set} $Y$ (or that $\cam$ is a \textbf{totally separating family }for $Y$) if  every $x, y\in Y$
are totally separated. 
\item We say that $\cam$ \textbf{well separates the set }$Y$
if for every $x\in Y$ and every $M\in \cam$ with $x\not \in M$ 
there exists a set $N\in \cam$ with $x\in N$ and $M\not \subseteq N$. 
\end{enumerate}
 \end{definition}
\begin{remark}
 The  notion of a separating family
is introduced by A. Renyi \cite{re} who has shown that the minimal size
of a separating family of a finite set $X$ is exactly $\left\lceil  \log_2|X|\right\rceil$.
\end{remark}
\begin{remark}
We note that for a family $\cam$ of sets and a set $Y\subseteq \bigcup\cam$
we have the obvious implications: \\
$\cam$ totally separates $Y$ $\Longrightarrow$ $\cam$ completely  separates
 $Y$ $\Longrightarrow$ $\cam$ well  separates $Y$
$\Longrightarrow$
$\cam$   separates $Y$.
\end{remark}

The main result of this section is the following:
\begin{theorem}\label{thlocset}
Let $X$ be a set and $\cam$ be a loc-lattice  of subsets of  $X$ which separates
the points of $X$.
\begin{enumerate}
 \item [(a)]  There exists a linear ordering $\leq$ of $X$ such that every $M\in \cam$ is
a convex subset of $X$.
 \item [(b)]   If moreover the family  $\cam$ well separates  $X$,  then for every
$A, B\in \cam$, $A\vee B=\mathrm{co}(A\cup B)$ and so $(\cam, \subseteq)$
is a sublattice of $(\mathrm{Co}(X),\subseteq)$.
\end{enumerate}

\end{theorem}

In order to prove Theorem \ref{thlocset} we shall state and prove several 
auxiliary lemmas. 
 
The following lemma summarizes the basic  properties of a loc-lattice of sets:

\begin{lemma}\label{lemp1}
Let $\cam$ be a loc-lattice of sets  and  $A, B, C\in \cam\setminus \{\emptyset\}$. Then
\begin{enumerate}
\item If the sets $A, B, C$ are pairwise incomparable then only one of the
relations $A\subseteq B\vee C$, $B\subseteq A\vee C$, $C\subseteq A\vee B$ occurs.\\
\item If $A, B, C$ are pairwise disjoint and $B\cap (A\vee C) \neq \emptyset$ then $B\subseteq A\vee C$.\\
\item If $A\cap C\neq \emptyset$ and $B\cap C\neq \emptyset$ then
$A\vee B\setminus A\cup B \subseteq C$.\\
\item If $C\cap (A\setminus B) \neq \emptyset$ and
$C\cap (B\setminus A)\neq \emptyset$ then
$A\cap B\subseteq C$.
\end{enumerate}
\end{lemma}
\begin{proof}
(1)\;\;  If $A\subseteq B\vee C$ and  $B\subseteq A\vee C$
then $A\subseteq (B\vee C) \cap (B\vee A)=B$, by Condition (5) of Definition \ref{def21}. \\
(2) \;\;If $A\subseteq B\vee C$ then by Condition (5) of Definition \ref{def21} we have that
$(A\vee B) \cap (A\vee C)=A$. But then $B\cap A = B\cap (A\vee C) \neq \emptyset$,
a contradiction. If $C\subseteq A\vee B$ then $C=(A\vee C)\cap (B\vee C)$
and so $B\cap(A\vee C)\subseteq (A\vee C)\cap (B\vee C)=C$, a contradiction.
Therefore by Condition (4) of Definition \ref{def21}  we must have $B\subseteq A\vee C$. \\
(3)\;\; Since  $A\cap C\neq \emptyset$ and $B\cap C\neq \emptyset$
we have that  $A\cup B\cup C= (A\vee B)\vee C\in \cam$. So,
$A\vee B\subseteq (A\cup B)\cup C$.\\
(4)\:\; We may suppose that $A\cap B\neq \emptyset$.
 The case $C\cap A\subseteq (C\cap B)\vee (A\cap B)$ is impossible since then
$C\cap A\subseteq B$ which contradicts the assumption $(C \setminus B)\cap A\neq \emptyset$.
The case $C\cap B\subseteq (C\cap A) \vee (A\cap B)$ is also impossible.
Therefore, $A\cap B\subseteq (C\cap A)\vee (C\cap B) \subseteq C$.   \end{proof}

\begin{definition}$ $
\begin{enumerate}
\item A\textbf{ ternary  relation} on a set $X$ is a subset $T$ of $X\times X\times X$. 
We shall use the notation
$(abc)_T$ instead of $(a, b, c)\in T$ and  $\neg(abc)_T$ instead of $(a, b, c)\not \in T$.
 \item If $\cam$ is a family of subsets of a set $X$ we define a ternary relation in $X$
setting $(abc)_\cam$ if and only if  for every $M\in \cam$ with $a, c\in M$ we have also  that
$b\in M$.  Let $$T_\cam=\{(a, b, c)\in X^3: (abc)_\cam\}.$$
 \item A binary relation $\precsim$ on a set $X$ is called a \textbf{quasi-order}
of it is reflexive and trasitive. Given a quasi- ordering
$\precsim$ on a set $X$ we set
$$T_\precsim=\{(a, b, c) \in X^3: a \precsim b \precsim c\;\; or \;\; c \precsim b \precsim a\}.$$
\item 
Let $\cam$ be a family of subsets of a set $X$ and $Y\subseteq X$.
We say that a quasi-order $\precsim$ of $Y$ of is $\cam$- \textbf{consistent}  
  if $T_\precsim \subseteq T_\cam$.
\end{enumerate}

\end{definition}

\begin{remark}

\noindent It is clear that in  order to prove Theorem \ref{thlocset} we must find an  $\cam$-consistent \textbf{linear} ordering
$\leq $ of $X$. 
\end{remark}

\begin{lemma}\label{lem04}
Let $\cam$ be a loc-lattice of subsets of a set  $X$. Then for every $a, b, c\in X$
either $(abc)_\cam$ or $(acb)_\cam$ or $(bac)_\cam$.\end{lemma}
\begin{proof} If $\neg(abc)_\cam$ and  $\neg(acb)_\cam$
then  there exist $M, N\in \cam$ such that $\{a, c\}\subseteq M$,
$\{a, b\}\subseteq N$, $b\not \in  M$ and $c\not \in  N$. Let $K\in \cam$ such that
$b, c\in K$ ( in particular, we can choose $K=M\vee N$).
 Then $K\cap (M\setminus N) \neq \emptyset$ and
$K\cap (N\setminus M) \neq \emptyset$. By Lemma \ref{lemp1} (4) we have that
$M\cap N\subseteq K$ and so $a\in K$.  \end{proof}

\begin{notation}
 If $A, B$ are subsets of a linear ordered set $(X, \leq)$ we write $A<B$ if for every
$a\in A$ and $b\in B$ we have that $a<b$. Similarly if $a\in X$ and $B\subseteq X$ we
write $a<B$ if $a<b$ for every $b\in B$.
\end{notation}

\begin{lemma}\label{lemlocsep} Let $\cam$ be a loc-lattice   of subsets of $X$
which  well-separates $X$
and $\leq$ an $\cam$-consistent linear ordering of $X$.
Then  $\cam$ is a sublattice of the lattice $\mathrm{Co}(X)$  of intervals of $(X, \leq)$.
\end{lemma}
\begin{proof} Clearly $\cam$ is a subset of $\mathrm{Co}(X)$.
 It remains to show that if $A, B\in \cam$ then $A\vee B=\mathrm{co}(A\cup B)$
where $\mathrm{co}(A)$ denotes the convex hull of a subset $A$ of $X$.
Since
$A\vee B$ is convex we have that
$\mathrm{co}(A\cup B)\subseteq A\vee B$. Suppose that
$\mathrm{co}(A\cup B)\neq A\vee B$. Then $A\cap B=\emptyset$ and therefore  either $A<B$ or
$B<A$.
Suppose that $A<B$.   Let
$c \in (A\vee B) \setminus \mathrm{co}(A\cup B)$. Then we must have that either
$c<A$ or $B<c$. Suppose that $B<c$ and let $C\in \cam$ with $c\in C$ and $C\not\subseteq B$.
The sets  $A, B, C$ are incomparable, $B\subseteq A\vee C$ and $C\subseteq A\vee B$ which contradicts Lemma \ref{lemp1} (1).
  \end{proof}

The following Theorem,
due to M. Altwegg \cite{al} (see also \cite{sh}),  is usefull
since it characterizes linear orderings of sets in terms of ternary relations.
Note that by Lemma \ref{lem04} $T_\cam$ always satisfies condition (3) of the Theorem of Altwegg.

\begin{theorem}[M. Altwegg]\label{thalt}
Let $T$ be a ternary relation in a set $X$  which satisfies the following postulates:
\begin{enumerate}
\item $(aba)_T$ if and only if $a=b$.\\
\item If $(abc)_T$ and $(bde)_T$ then either $(cbd)_T$ or $(eba)_T$.\\
\item  For every $a, b, c\in X$ either $(abc)_T$ or $(bca)_T$
or $(cab)_T$.
\end{enumerate} Then there exists a
  linear ordering
$\leq$ of $X$ such that for every $x, y, z\in X$ we have that
$(xyz)_T$ if and only if $(xyz)_{T_\leq}$. Moreover, every other linear ordering $\leq'$
of $X$ satisfying the preceding condition is equal to $\leq$ or to the inverse order $\leq^*$
of $\leq$.
\end{theorem}

The ternary relation $T_\cam$   does not satisfy in general
the conditions of  Theorem \ref{thalt} since
for a loc-lattice  $\cam$  of subsets of $X$
there are probably  many and even non-isomorphic $\cam$-consistent linear orderings.

 Indeed, the
simplest example of a loc-lattice   of subsets of a set $X$   is a  chain of subsets of  $X$.
Suppose that $X=\mathbb{N}$ and that
$\cam=\{A_n: n\in \mathbb{N}\}$ where $A_n=\{1, \dots, n\}.$
Then the usual ordering $1\leq 2\leq 3\leq \dots $ of $\mathbb{N}$ and
the ordering $\leq'$ given by
$$\dots\leq'2n+2 \leq' 2n\leq'\dots4\leq' 2\leq' 1 \leq' 3 \leq'\dots \leq'2n-1 \leq'2n+1\leq'\dots $$
are non-isomorphic $\cam$-consistent   linear orderings of $\mathbb{N}$.

\begin{definition}
 Let $X$ be a set and $\cam$ a family of subsets of $X$.
If $(x_i)_{i=1}^n$ is a sequence of pairwise distinct points of $X$
we say that $(M_i)_{i=1}^n$ is a \textbf{representative family} for
$(x_i)_{i=1}^n$ if for every $i, j\in \{1, \dots, n\}$ with  $j\neq i$ we have that
 $M_i\in \cam$, $x_i\in M_i$
and $x_j\not \in M_i$.
\end{definition}

It is clear that if the family
$\cam$ is closed under finite intersections then for every finite sequence
$(x_i)_{i=1}^n$ of pairwise completely separated points of $X$ there exists
a  \textit{representative family} for
$(x_i)_{i=1}^n$.

\begin{lemma}\label{lemp3} Let $\cam$ be a loc  lattice of subsets of a set $X$, $a, b, c$ three    distinct points of $X$
which  are pairwise completely separated. Let $(A, B, C)\in \cam^3$
a representative family for the
triple $(a, b, c)$.
Then
\begin{enumerate}
\item [$(\mathrm{i})$]  $(abc)_\cam$ holds if and only if $B\subseteq A\vee C$.\\
\item [$(\mathrm{ii})$] If $B\subseteq A\vee C$ then
$A\cap C=\emptyset$, $a\not \in B\vee C$ and $c\not \in A\vee B$.
\end{enumerate}
 \end{lemma}
\begin{proof}
$(\mathrm{i})$:\;\;\;\;We shall show that if  $(abc)_\cam$ holds then
for every representative family $(A, B, C)$ of $(a, b, c)$ we have that $B\subseteq A\vee C$.
If $A\subseteq B\vee C$ then $A=(A\vee B)\cap (A\vee C)$ and so
$A\cap B= B \cap (A\vee C)$. Since $(abc)_\cam$ we have that $b\in A\vee C$ and so
$b\in A$, a contradiction. In the same manner the case $C\subseteq A\vee B$
is excluded. Therefore, we must have $B\subseteq A\vee C$.
Suppose now that there exists a  representative family  $(A, B, C)$ for $(a, b, c)$ such that
$B\subseteq A\vee C$ and $(abc)_\cam$ does not hold. Then by Lemma \ref{lem04}
either $(bac)_\cam$ or $(bca)_\cam$ and by (a), either $A\subseteq B\vee C$
or $C\subseteq A\vee B$, but this contradicts Lemma \ref{lemp1} (1).

\medskip

\no$(\mathrm{ii})$:\;\;\;\; Suppose that  $A\cap C\neq \emptyset$ then $A\vee C=A\cup C$ and
since $B\subseteq A\vee C$ we have that either $b\in A$ or $b\in C$,
which contradicts the assumption that $(A, B, C)$
is a representative family for $(a, b, c)$.   If $a\in B\vee C$ since $B=(A\vee B)\cap (B\vee C)$
we shall have $a\in B$, a contradiction. Similarly, we see that  $c\not \in A$.  \end{proof}

\begin{lemma}\label{lem02}
Let $\cam$ be a loc-lattice of subsets of a set $X$ which completely separates
$X$. Then $\cam$ totally separates $X$.
\end{lemma}
\begin{proof} Let $a, b$ be distinct points of $X$ and $(A_1, B_1)$
a representative family for $(a, b)$. Suppose that $A_1\cap B_1 \neq \emptyset$
and let $c\in A_1\cap B_1$. Then we may choose a representative family
$(A, B, C)$ for $(a, b, c)$ with $A\subseteq A_1$, $B\subseteq B_1$, $C\subseteq A_1\cap B_1$.
Since $C\vee A\subseteq C\vee A_1=A_1$ we cannot have $B\subseteq C\vee A$. Similarly,
the case $A \subseteq C\vee B$ is impossible. Therefore by Property (3) of Definition \ref{def21} we must have
that $C\subseteq A\vee B$.  By Lemma \ref{lemp3} $(\mathrm{ii})$  we conclude that
$A\cap B=\emptyset$. \end{proof}

\begin{lemma}\label{lemp4}
Let $\cam$ be a loc  lattice of subsets of a set $X$, $a, b, c, d, e$ five  distinct points of $X$
such that every two of them are completely separated and
  $(A, B, C, D, E)$
 a representative family for $(a, b, c, d, e)$. Then\\
\begin{enumerate}
\item $B\subseteq A\vee C$ and $C\subseteq B\vee D$ then
$B\subseteq A\vee D$ and $C\subseteq A\vee D$.\\
\item $B\subseteq A\vee C$ and $D\subseteq B\vee C$ then 
$B\subseteq A\vee D$ and $D\subseteq A\vee C$.\\
\item If $B\subseteq A\vee C$ and $D\subseteq B\vee E$ then
either $B\subseteq C\vee D$ or $B\subseteq A\vee E$.\\
\end{enumerate}
\end{lemma}
\begin{proof} (1) \;\; Suppose that $B\subseteq A\vee C$ and $C\subseteq B\vee D$. 
Since $B\subseteq A\vee C$ and $A, B, C$ are pairwise incomparable we
have by Property (5) of Definition \ref{def21}  that $B=(A\vee B)\cap (B\vee C)$  and so  $c\not \in A\vee B$. 

Clearly, $A\vee B\subseteq(A\vee D) \vee (B\vee D)=(A\vee D) \cup (B\vee D)$.
 If $D\subseteq A\vee B$ then  $(A\vee D) \cup (B\vee D)\subseteq A\vee B$,
and so $A\vee B=(A\vee D) \cup (B\vee D)$.
But  $C\subseteq B\vee D$ which implies that 
$C\subseteq A\vee B$, which contradicts the fact that $c\not\in A\vee B$. 
So, the case $D\subseteq A\vee B$ is impossible. 

Also, $A\vee C\subseteq (A\vee D)\cup(D\vee C)$ and since $B\subseteq
A\vee C$ then  $b\in A\vee D$ or $b\in D\vee C$. 
Since $C\subseteq B\vee D$ by Lemma \ref{lemp3} $(\mathrm{ii})$ we have that $b\not\in C\vee D$ and so $ b\in A\vee D$. 
If it was true that $A\subseteq B\vee D$ then, 
again by Lemma \ref{lemp3} $(\mathrm{ii})$, $b\not\in A\vee D$, a contradiction.
Therefore, the case $A\subseteq B\vee D$ is also impossible.
By Definition \ref{def21},
either $A\subseteq B\vee D$ or $B\subseteq A\vee D$ or $D\subseteq A\vee B$. So, $B\subseteq A\vee D$.\\(2)\;\;Similar to (1).\\ 
(3)\;\;\;Let $B\subseteq A\vee C$ and $D\subseteq B\vee E$. 
Suppose that  $B\not\subseteq C\vee D$. Then either $C\subseteq B\vee D$
or $D\subseteq B\vee C$. \\ 
 \textbf{Case (3a): }  $C\subseteq B\vee D$.  Since,  $B\subseteq A\vee C$
and $C\subseteq B\vee D$
we have by (1) that $B\subseteq  A\vee D$. Since $B\subseteq A\vee D$ and 
$D\subseteq B\vee E$ we have again by (1) that  $B\subseteq A\vee E$.\;\; \\
 \textbf{Case (3b): } $D\subseteq B\vee C$. Since $B\subseteq A\vee C$
and $D\subseteq B\vee C$ we have by (2) that $B\subseteq A\vee D$.
Since $B\subseteq A\vee D$ and $D\subseteq B\vee E$ we have by (1)
that $B\subseteq A\vee E$.
\end{proof}

\begin{lemma}\label{prop}
Let $\cam$ be a loc-lattice of subsets of a set $X$
 which completely separates the points of a subset $Y$ of $X$. 
Then there exists a  consistent linear ordering $\leq$ of $Y$.
Moreover, every other consistent linear ordering $\leq'$ of $Y$
  is equal to $\leq$ or equal to the inverse order of $\leq$.
\end{lemma}
\begin{proof}
Let $\cam_Y=\{M\cap Y: M\in \cam\}$.  By Lemmas \ref{lem04},  \ref{lemp3} and \ref{lemp4}
we see  that the ternary relation  $T_{\cam_Y}$ of $Y$ satisfies the conditions of Altwegg's Theorem \ref{thalt}.
So we may find a $\cam_Y$-consistent linear ordering of $Y$ which will be
$\cam$-consistent too.
\end{proof}

\no By Lemmas \ref{prop} and   \ref{lemlocsep} we have as an immediate corollary the following:

\begin{theorem}\label{thsep}
Let $X$ be a set and $\cam$ be a loc-lattice  of subsets of a set  $X$ which completely separates
$X$. Then there exists a linear ordering $\leq$ of $X$ such that
 $\cam$ is a sublattice of the lattice $(\mathrm{Co}(X), \subseteq)$ of all intervals
of $(X, \leq)$. Moreover, every other linear ordering $\leq'$
of $X$ satisfying the preceding condition is equal to $\leq$ or to the inverse order $\leq^*$
of $\leq$.
\end{theorem}

The general case of Theorem \ref{thlocset}, when the family $\cam$ does not completely
separate the set $X$ is much more complicated.

If $X, Y$ are sets and $R\subseteq X\times X$,
$S\subseteq Y\times Y$ are binary relations we say that 
$(X, R)$ is an \textit{extension} of $(Y, S)$ if $X\supseteq Y$
and $R\supseteq S$.

In order to prove the theorem 
we shall define for every ordinal $\alpha$ a subset 
$Y_\alpha$ of $X$, a linear ordering $\leq_\alpha$ of $Y_\alpha$ 
and a quasi-ordering $\precsim_\alpha$ of $X$ such that the following conditions are 
fullfilled:\\

\begin{itemize}\label{aaa}
\item[(a)] For every ordinal $\alpha$, $\leq_\alpha$ is an $\cam$-compatible
linear ordering of $Y_\alpha$ and $\precsim_\alpha$ is an $\cam$-compatible
quasi-ordering of $X$.\\
\item [(b)] For every ordinal $\alpha$, $(X, \precsim_\alpha)$ is an extension 
of $(Y_\alpha, \leq_\alpha)$.\\
 \item [(c)] If $\alpha, \beta$ are ordinals and  $\alpha\leq \beta$ then
$Y_\alpha\subseteq Y_\beta$. \\
\item[(d)]  If $Y_\alpha\neq X$ then $Y_{\alpha+1}\neq Y_\alpha$. 
\end{itemize}
\bigskip

If a such  construction of $(Y_\alpha, \leq_\alpha, \precsim_\alpha)$ 
is possible for every ordinal $\alpha$ then there 
exists an ordinal $\alpha$ such that $Y_\alpha=X$ and so 
$\leq_\alpha=\precsim_\alpha$  will be an $\cam$-compatible
linear ordering of $X$ which proves Theorem \ref{thlocset}.

\noindent \underline{\textbf{Step} 0}.

\medskip

\no We start by some auxiliary definitions:
\begin{definition}
 A \textbf{section}  of a linearly ordered set  $(X, \leq) $ 
is a  pair $\mathbf{S}=(A, B)$ of subsets of $X$ such that $A\cup B=X$ and $A<B$.  
The  pairs $(\emptyset, X)$
and $(X, \emptyset)$ are considered as sections too.
The set of all sections of a linearly ordered set
$(X, \leq)$ is denoted by  $\mathfrak{S}(X, \leq)$.
 Given a section $\mathbf{S}=(A, B)$  we also denote by 
$\mathbf{S}^{(1)}$ the first member $A$ of the section
and by $\mathbf{S}^{(2)}$  the second member $B$ of the section. 
\end{definition}
\begin{definition}
We  define a binary relation $\mathbf{L}\subseteq X\times X$ on $X$ so that for $x, y\in X$ $x\mathbf{L} y$ holds  if and only 
if for   every  $M\in \cam$ such that  $x\in M$ we have that $y\in M$. 
We write $\neg(x\mathbf{L} y)$ if $x\mathbf{L} y$ does not hold.
Two points $x, y$ of $X$ are said to be \textbf{independent} if $\neg (x\mathbf{L} y)$ and $\neg(y\mathbf{L} x)$. A subset $Y$ of $X$ is said to be \textit{independent} if every two distinct elements of $Y$ are independent, that is if $Y$ is  completely separated 
by the family $\cam$.
\end{definition}
\begin{remark}\label{rem11}
It is trivial that  $\mathbf{L}$ is reflexive and  transitive, and since
  $\cam$ separates the set $X$ it is also antisymmetric;
 that is $\mathbf{L}$ is  a partial order on  $X$.
Note also that if $Y$ is an independent subset of $X$ then by Proposition
\ref{prop} there exists an $\cam$-consistent  linear ordering of the set $Y$.
\end{remark}

If $\mathcal{A}$ is a family of independent subsets of $X$, 
linearly ordered by $\subseteq$, then $\bigcup\mathcal{A}$ is an independent 
set. Therefore,  by Zorn's Lemma the set $X$ contains a maximal independent subset $Y$. This set admits, by Lemma \ref{prop} an $\cam$-consistent linear ordering $\leq$. 

We fix a maximal independent subset  $Y=Y_0$  of $X$ and $\leq=\leq_0$ an $\cam$-consistent linear ordering of $Y$.

If $|Y|=1$ then  if $x, y\in X$ either $x\mathbf{L}y$ or $y\mathbf{L}x$, otherwise 
the set $\{x, y\}$ would be independent, which contradicts our assumption that 
a maximal independent set has size equal to 1. So, by Remark \ref{rem11}, the relation $\mathbf{L}$ is an $\cam$-consistent linear ordering of $X$ and Theorem \ref{thlocset} has been  proved.  So, we may assume that $|Y|\geq 2$.

\begin{definition}
 For every $x\in X\setminus Y$  we set
$$
L(x)=L_Y(x)=\{y\in Y: x\mathbf{L} y\}, $$
$$ M(x)=M_Y(x)=\{y\in Y: y\mathbf{L} x\}.$$
\end{definition}

Note that since $Y$ is a maximal independent subset of $X$,  for every $x\in X\setminus Y$ we have that $L(x)=\emptyset$ if and only if $M(x)\neq \emptyset$ and
that $L(x)$ is a convex subset of $(Y, \leq)$. 

\begin{lemma}\label{lemM}
For every $x\in X$   we have that $|M(x)|\leq 2$ and if
$M(x)=\{y, z\}$ with $y<z$ then $z$ is the immediate successor of $y$.
\end{lemma}
\begin{proof}
Supose that there exist three distinct elements $a, b, c\in M(x)$.
Since $a, b, c\in Y$ they are independent, and so there 
exists a representative family $(A, B, C)$ for $(a, b, c)$.
We may assume that  $A\subseteq B\vee C$ and then by  Lemma \ref{lemp3} (ii) 
we must have that $B\cap C=\emptyset$. But $x\in B\cap C$, a contradiction. 

If there exists $w\in Y$ with $y<w<z$ then again by  Lemma \ref{lemp3}  
we can find a representative family  $(A, B, C)$ for $(y, w, z)$ such that 
$A\cap C=\emptyset$, which contradicts the fact that $x\in A\cap C$.
\end{proof}

\begin{notation}
 Let $(Y, \leq)$ be  a linearly ordered set. If $A\subseteq Y$ we set
$$A^-=\{y\in Y:\; \text{for every}\; a\in A, \;y<a\},$$
$$ A^+=\{y\in Y:\; \text{for every}\; a\in A,\; a<y\}.$$
\end{notation}

Using Lemma \ref{lemM} we may  classify all the points of $X\setminus Y$ in eight 
mutually exclusive classes or types by the following rules:

\begin{definition}
 A point $x\in X\setminus Y$ is said to be a point  of
\\[1ex]
$\mathbf{Type\;  I_a}$\;\; if $L(x)\neq \emptyset$, $L^+(x)\neq \emptyset$ and $L^-(x)\neq \emptyset$,\\[1ex]
$\mathbf{Type\; I_b}$\;\; if $L(x)\neq \emptyset$, $L^+(x)\neq \emptyset$ and $L^-(x)= \emptyset$,\\[1ex]
$\mathbf{Type \;I_c}$\;\; if $L(x)\neq \emptyset$, $L^+(x)= \emptyset$ and $L^-(x)\neq \emptyset$,\\[1ex]
$\mathbf{Type \;I_d}$ \;\;if $L(x)=Y$,\\[1ex]
$\mathbf{Type \;II_a}$\;\; if $|M(x)|=2$,\\[1ex]
$\mathbf{Type \;II_b}$ \;\;if $M(x)=\{y\}$ and $y$ is neither the first nor the  last element of $Y$, \\[1ex]
$\mathbf{Type \;II_c}$ \;\;if $M(x)=\{y\}$ and $y$ is  the first  element of $Y$,\\[1ex]
$\mathbf{Type \;II_d}$\;\; if $M(x)=\{y\}$ and $y$ is  the  last element of $Y$.\\[1ex]
\end{definition}

\begin{definition}
 For every  $x\in X\setminus Y$  with $L(x)\neq \emptyset$. We set
$$\begin{array}{llll}
L_1(x)=\{y\in L(x):  \text{there exists a} \;\;y' \in L^-(x)
\; \text{with}\; \neg(yxy')_\cam\}, \\[1.5 ex]
 L_2(x)=\{y\in L(x): \text{there exists a }  \; y'\in L^+(x) \; \text{with}\; \neg(yxy')_\cam\}.\\
\end{array}
$$
\end{definition}

Note that the linearly ordered set under consideration is $(Y, \leq)$, so the sets   $L^-(x), L^+(x)$ are always subsets of $Y$. 

Our aim in this stage of the proof is to associate to  every $x\in X\setminus Y$  a section $\mathbf{S}_x=(A_x, B_x)$ of $Y$. 
In order to do this for the points  of type $\mathrm{I_a}$ we  first need a lemma. 
\medskip
\begin{lemma}\label{lemtype1}
Let $x\in X\setminus Y$  such that $L(x)\neq \emptyset$ and $L(x)\neq Y$. Then 
\begin{enumerate}
 \item The set
$L_1(x)$ is an initial segment of $L(x)$ (if it is not empty) and $L_2(x)$ is a final segment of $L(x)$ (if it is not empty). 
\item
Moreover, 
if $L^+(x)\neq \emptyset$ and $L^-(x)\neq \emptyset$ then 
$$L_1(x)\cap L_2(x)=\emptyset \;\; \text{and} \;\; L_1(x)\cup L_2(x)=L(x). $$
\end{enumerate}
\end{lemma}
\begin{proof}
Let $y\in L_1(x)$ and $y'\in L(x)$ with $y'<y$. Then there exists a $z<L(x)$ and $M\in \cam$
with $z, y\in M$ and $x\not \in M$. Since $z<y'<y$ and $\leq$ is consistent we have that
$y'\in M$ and so $y'\in L_1(x)$. This shows that $L_1(x)$ is an initial segment of $L(x)$. By the same reasoning we
show that $L_2(x)$ is a final segment of $L(x)$.

Suppose that there exists $y\in L_1(x)\cap L_2(x)$. Then, by the definition of the sets  $L_1(x)$
and $ L_2(x)$,  there exist  $y_1, y_2\in Y$
with $y_1<L(x)<y_2$ so that $\neg(yxy_1)_\cam$ and $\neg(yxy_2)_\cam$.
Since $\neg(yxy_1)_\cam$ and $\neg(yxy_2)_\cam$
there exist $N_y, N_{y_1}, N_{y_2}\in \cam$ such that 
$y\in N_{y}, y_1\in N_{y_1}, y_2\in N_{y_2}$ such
that $x\not\in   N_{y}\vee N_{y_1}$ and $x\not\in   N_{y}\vee N_{y_2}$.
Since $y_1, y_2,y\in Y$ these points are independent and so 
we may also assume that $(N_y, N_{y_1}, N_{y_2})$
is a representative family for $(y, y_1, y_2)$.
Since $y_1, y_2\not\in L(x)$ there exist $N^{(1)}_x, N^{(2)}_x\in \cam$ such that 
$x\in N_x$ and $y_1\not \in N^{(1)}_x, y_2\not \in N^{(2)}_x$.

We set $M_1= N_{y}\vee N_{y_1}$, $M_2=N_{y}\vee N_{y_2}$, $M=N^{(1)}_x\cap N^{(2)}_x$.
Since $y_1<y<y_2$ we have by Lemma \ref{lemp3} $(\mathrm{i})$ that 
$y_2\not \in M_1$ and $y_1\not\in M_2$ which means that   $(M_1, M_2, M)$ is a representative triple for $(y_1, y_2, x)$ such that   $y_1, y\in M_1$, $y_2, y\in M_2$. Since $\neg(yxy_1)_\cam$ and $\neg(yxy_2)_\cam$,
by  Lemma \ref{lemp3} $(\mathrm{i})$, we have that  $M_1\not \subseteq M\vee M_2=
M\cup M_2$ and  $M_2\not \subseteq M\vee M_1=
M\cup M_1$.   But then we must have that  $M\subseteq M_1\vee M_2$
and so by Lemma \ref{lemp3} $(\mathrm{ii})$.
 $M_1\cap M_2= \emptyset$. But $y\in M_1\cap M_2$, a contradiction.  

Finally suppose that there exists a $y\in L(x)\setminus (L_1(x)\cup L_2(x))$.
Since $L^+(x)\neq \emptyset$ and $L^-(x)\neq \emptyset$, we may find  $y_1, y_2\in Y$
with $y_1<L(x)<y_2$ and let  $(M_1, M_2, M)$ be  a representative triple for $(y_1, y_2, y)$.
Since $y\in L(x)$ we have that $x\mathbf{L}y$ and since the family 
separates $X$ we have that  $\neg y\mathbf{L}x$, that is we may 
choose the set  $M$ such that  $x\not \in M$. Since $y\not \in L_1(x)$ we have that $x\in M\vee M_1$ and
since $y\not\in L_2(x)$ we have that $x\in M\vee M_2$. But then
  $x\in (M\vee M_1)\cap (M\vee M_2)=M$,
 a contradiction.
\end{proof}

The precending lemma shows that if   $x\in X\setminus Y$ 
 is a point of  type $\mathrm{I_a}$ then $(L_1(x), L_2(x)$ is a 
section of $(L(x), \leq)$ which can be extended to a section $\mathbf{S}_x=(A_x, B_x)$ of 
$(Y, \leq)$. In the following definition we associate a section $(A_x, B_x)$
to every point $x$ which is of type $\mathrm{I_a}$or $\mathrm{I_b}$ or $\mathrm{I_c}$,

\begin{definition}[\textbf{Sections associated to points of type $\mathrm{I_a}$}]\label{defIa}\hfill
\begin{itemize}
 \item [(i)] If $x$ is a point of type $\mathrm{I_a}$
such that $L_1(x)\neq\emptyset $ and  $L_2(x)\neq\emptyset $ we set
 $$\begin{array}{llll}
A_x=L_1^-(x)\cup L_1(x)=L_2^-(x),\\ [1ex]
B_x=L_1^+(x)=L_2(x)\cup L_2^+(x).
\end{array}$$
\begin{center}
\includegraphics[scale=0.19]{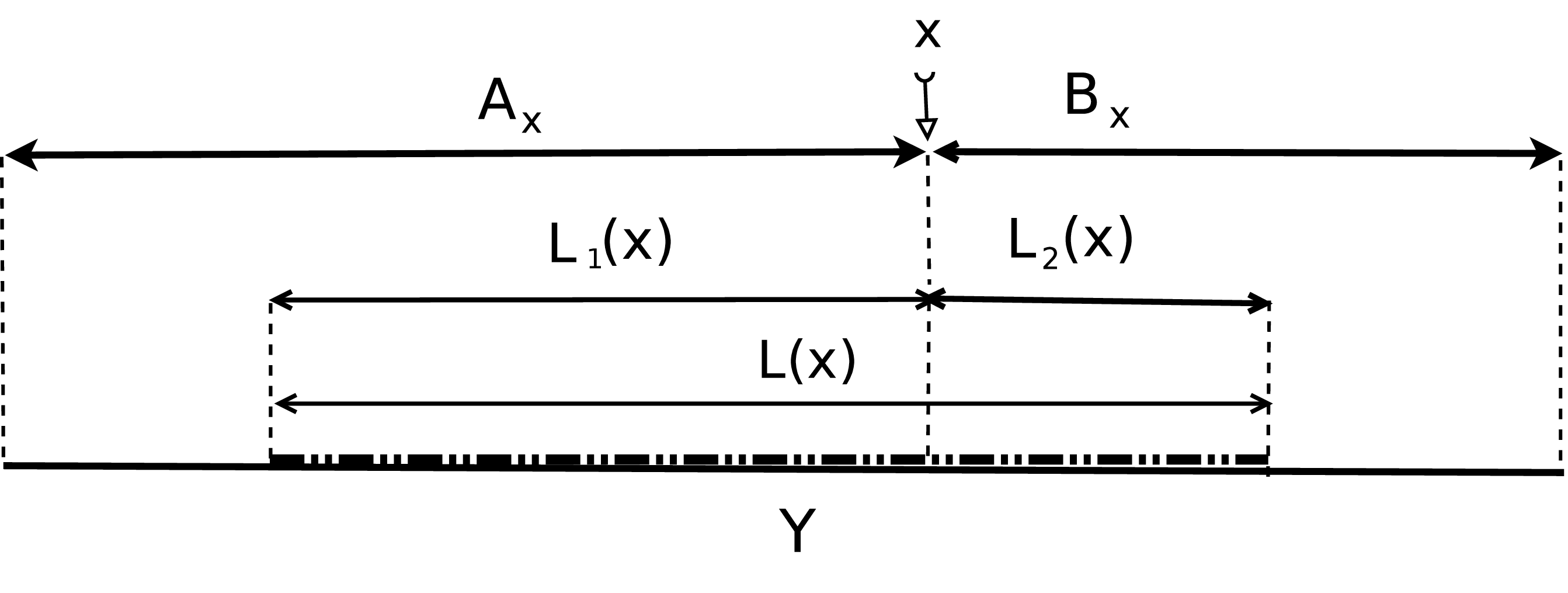}\\
\textbf{Fig.2}\;\; \begin{small}Points of type $I_a$ 
such that $L_1(x)\neq\emptyset $ and  $L_2(x)\neq\emptyset $ .                                                              \end{small}\end{center}
\item [(ii)] If $x$ is a point of type $\mathrm{I_a}$
such that $L_1(x)=\emptyset $, then by Lemma \ref{lemtype1}
we have that 
 $L_2(x)=L(x)$ 
and in this case we set
 $$\begin{array}{llll}
A_x=L^-(x), \;\;
B_x=L(x)\cup L^+(x).
\end{array}$$
\begin{center}
\includegraphics[scale=0.19]{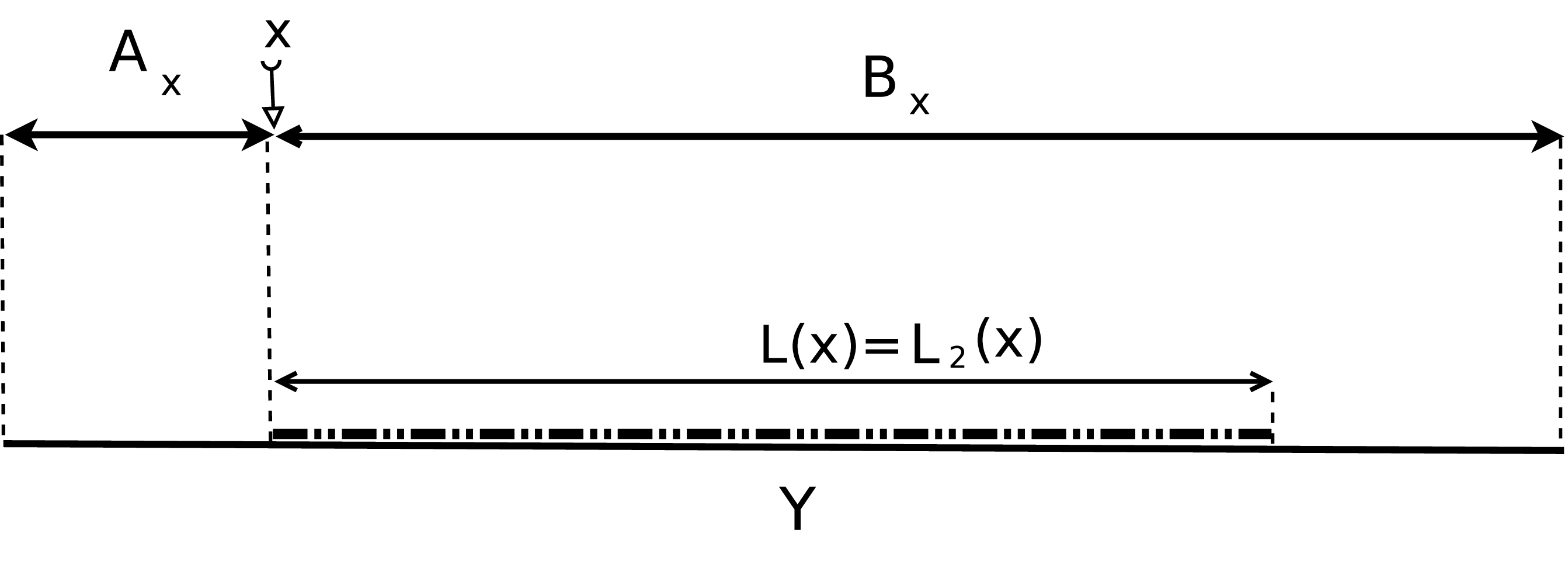}\\
\textbf{Fig.3}\;\; \begin{small}Points of type $I_a$ 
such that such that $L_1(x)=\emptyset $ .                                         \end{small}\end{center}
\item [(iii)] If $x$ is a point of type $\mathrm{I_a}$
such that $L_2(x)=\emptyset $, then by Lemma \ref{lemtype1}
we have that 
 $L_1(x)=L(x)$ 
and in this case we set
 $$\begin{array}{llll}
A_x=L^-(x)\cup L(x), \;\;
B_x=L^+(x).
\end{array}$$
\begin{center}
\includegraphics[scale=0.19]{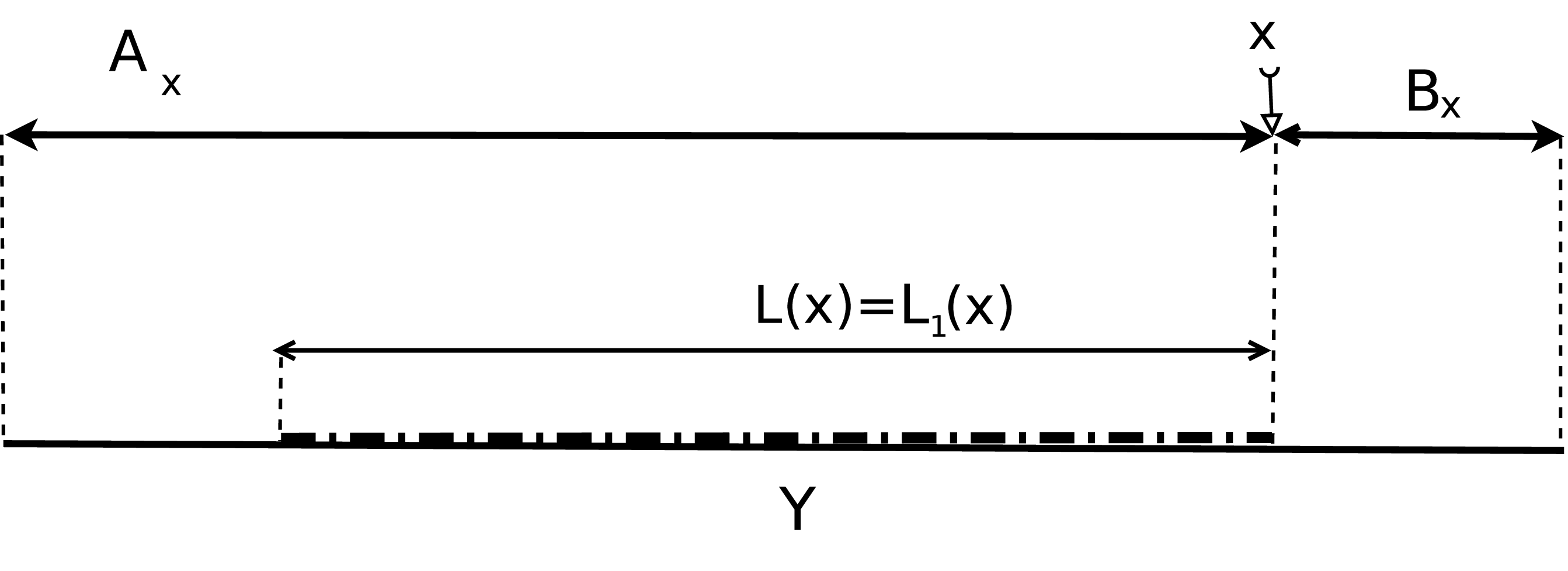}\\
\textbf{Fig.4}\;\; \begin{small}Points of type $I_a$ 
such that $L_2(x)=\emptyset $ .                               \end{small}\end{center}
\end{itemize}
\end{definition}

\medskip
 \begin{definition}[\textbf{Sections associated to points of type $\mathrm{I_b}$}]\hfill

If $x$ is a point of type $\mathrm{I_b}$, that is 
$L(x)$ is a proper initial segment of $Y$,  then we have that  $L_1(x)=\emptyset$.

\begin{itemize}
 \item [(i)] If  $L_2(x)\neq \emptyset$ and $L_2(x)\neq L(x)$ we set
 $$\begin{array}{llll}
A_x=L_2^-(x),\;\;\;
B_x=L_2(x)\cup L_2^+(x).
\end{array}$$
 \begin{center}
\includegraphics[scale=0.19]{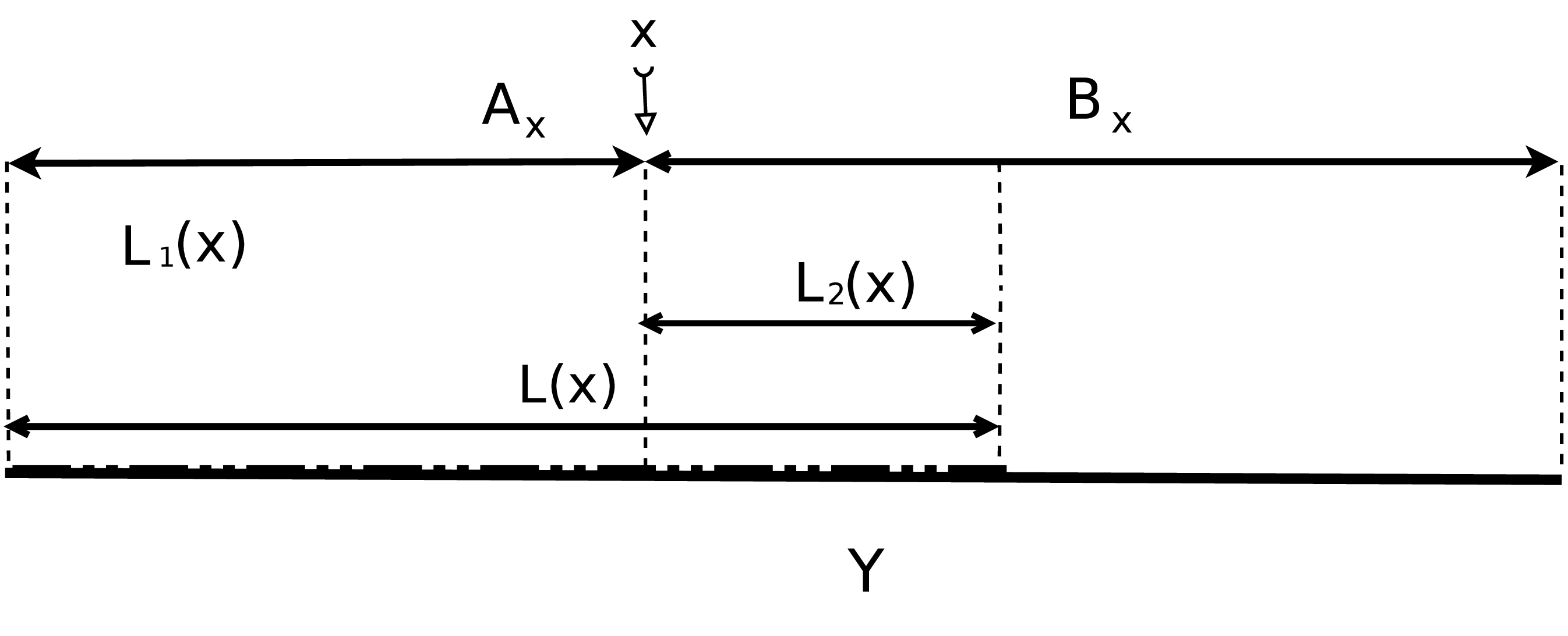}\\
\textbf{Fig.5}.\;\; \begin{small}Points of type $I_b$ such that $L_2(x)\neq \emptyset$ and $L_2(x)\neq L(x)$                                                                                                 \end{small}.
\end{center}
\medskip
\item [(ii)] If  $L_2(x)= L(x)$ we set
$$\begin{array}{llll}
A_x=\emptyset,\;\;\;
B_x=Y.
\end{array}$$ 
\begin{center}
\includegraphics[scale=0.2]{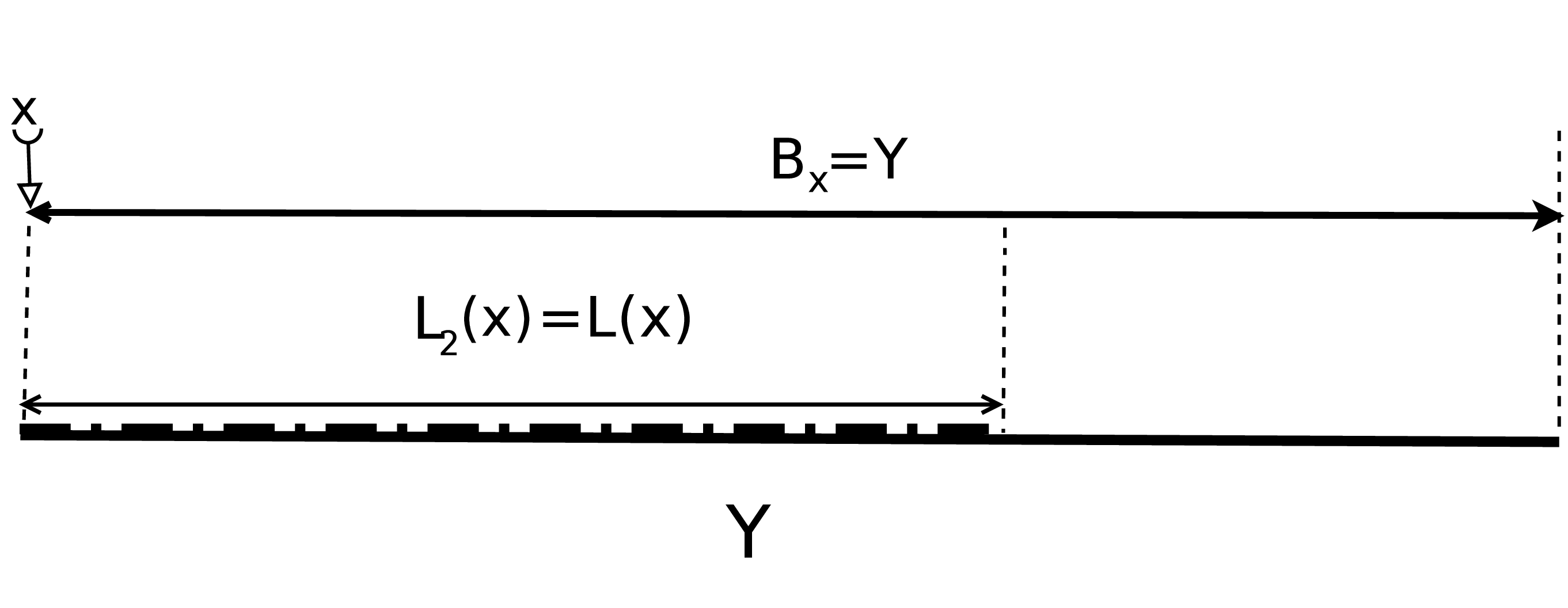}\\
\textbf{Fig.6}.\;\; \begin{small}Points of type $I_b$ such that  $L_2(x)=L(x)$ .                                                                   \end{small}
\end{center}
\medskip
 \item [(iii)] If $L_2(x)=\emptyset$, that is for every 
$y\in L(x)$ and every $y'\in Y$ with $y'>L(x)$ we have that 
$(yxy')_\cam$, then we set
$$\begin{array}{llll}
A_x=L(x),\;\;\;
B_x=L^+(x).
\end{array}$$ 

\begin{center}
\includegraphics[scale=0.19]{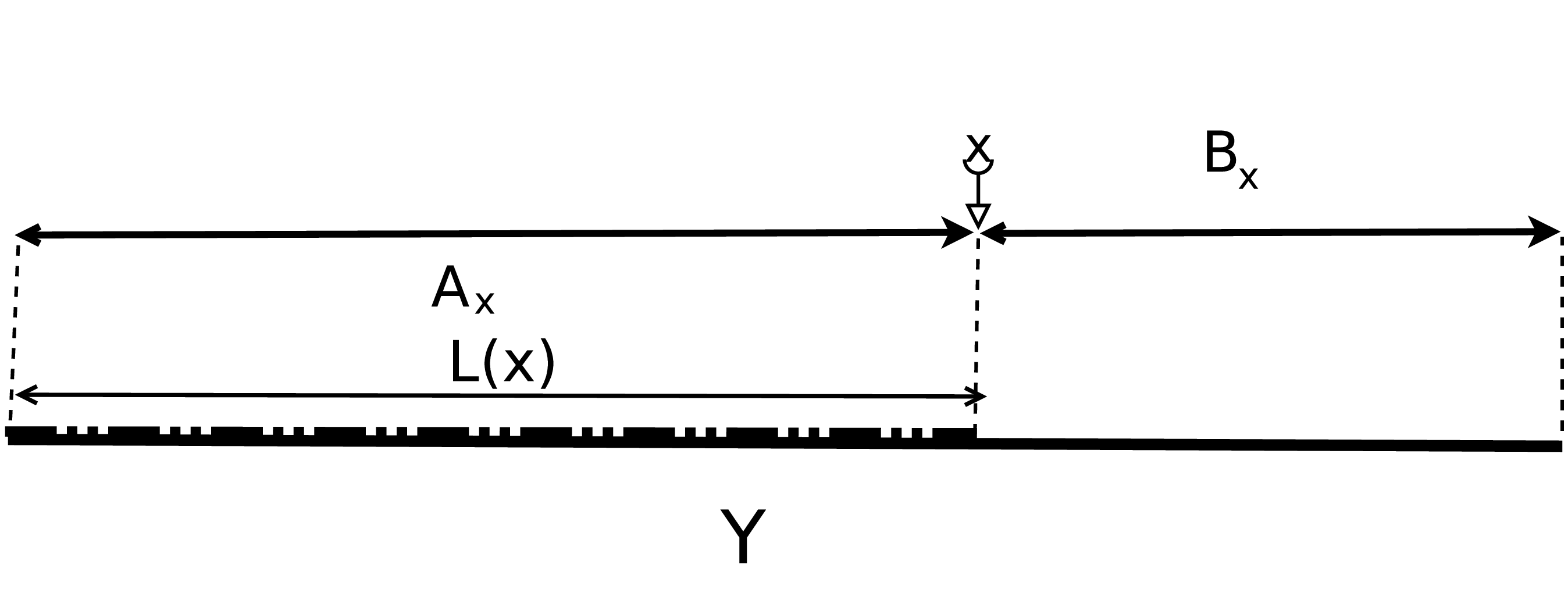}\\
\textbf{Fig.7}.\;\; \begin{small}Points of type $I_b$ such that  $L_2(x)=\emptyset$                                                                        \end{small}.
\end{center}
\medskip
\end{itemize}
\end{definition}
\medskip

\begin{definition}[\textbf{Sections associated to points of type $\mathrm{I_c}$}]\hfill

If $x$ is a point of type $\mathrm{I_c}$, that is 
$L(x)$ is a proper final segment of $Y$,  then we have that  $L_2(x)=\emptyset$.

\begin{itemize}
 \item [(i)] If  $L_1(x)\neq \emptyset$ and $L_1(x)\neq L(x)$ we set
 $$\begin{array}{llll}
A_x=L_1^-(x),\;\;\;
B_x=L_1(x)\cup L_2^+(x).
\end{array}$$
 \begin{center}
\includegraphics[scale=0.19]{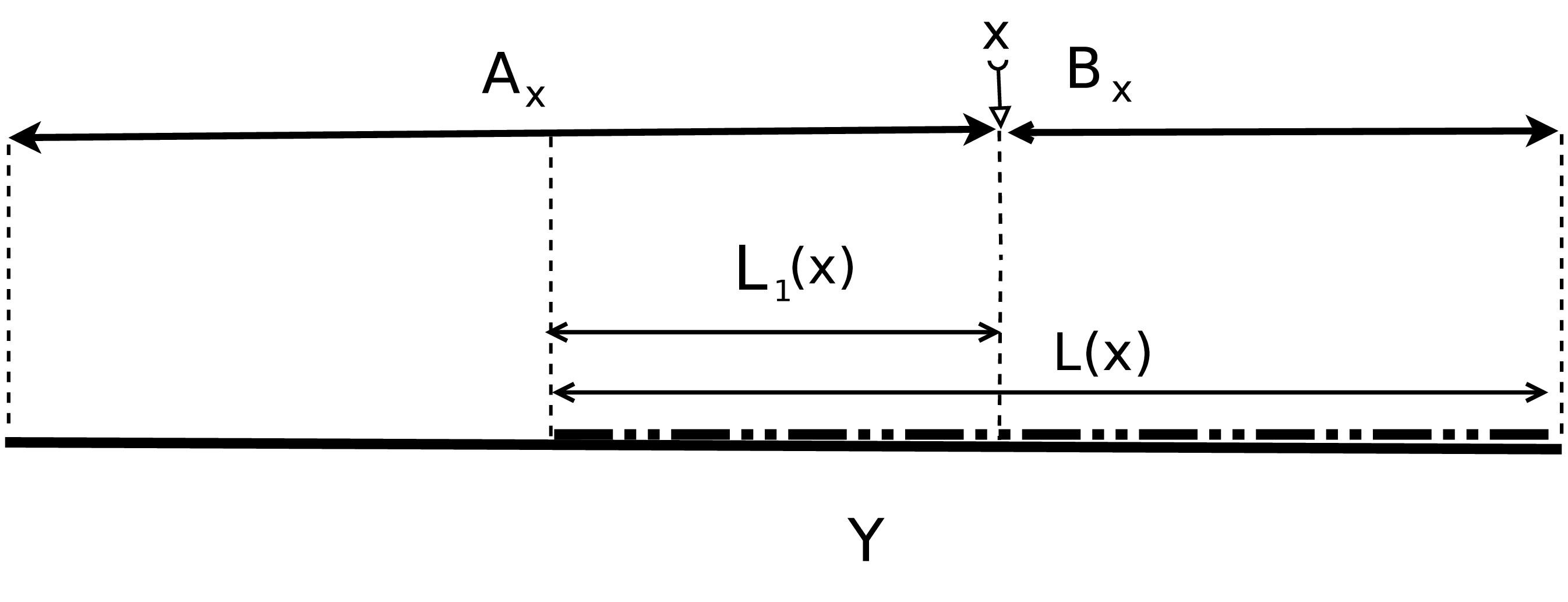}\\
\textbf{Fig.8}.\;\; \begin{small}Points of type $I_c$ such that $L_1(x)\neq \emptyset$ and $L_2(x)\neq L(x)$                                                                                                 \end{small}.
\end{center}
\medskip
\item [(ii)] If  $L_1(x)= L(x)$ we set
$$\begin{array}{llll}
A_x=Y,\;\;\; B_x=\emptyset.
\end{array}$$ 
\begin{center}
\includegraphics[scale=0.2]{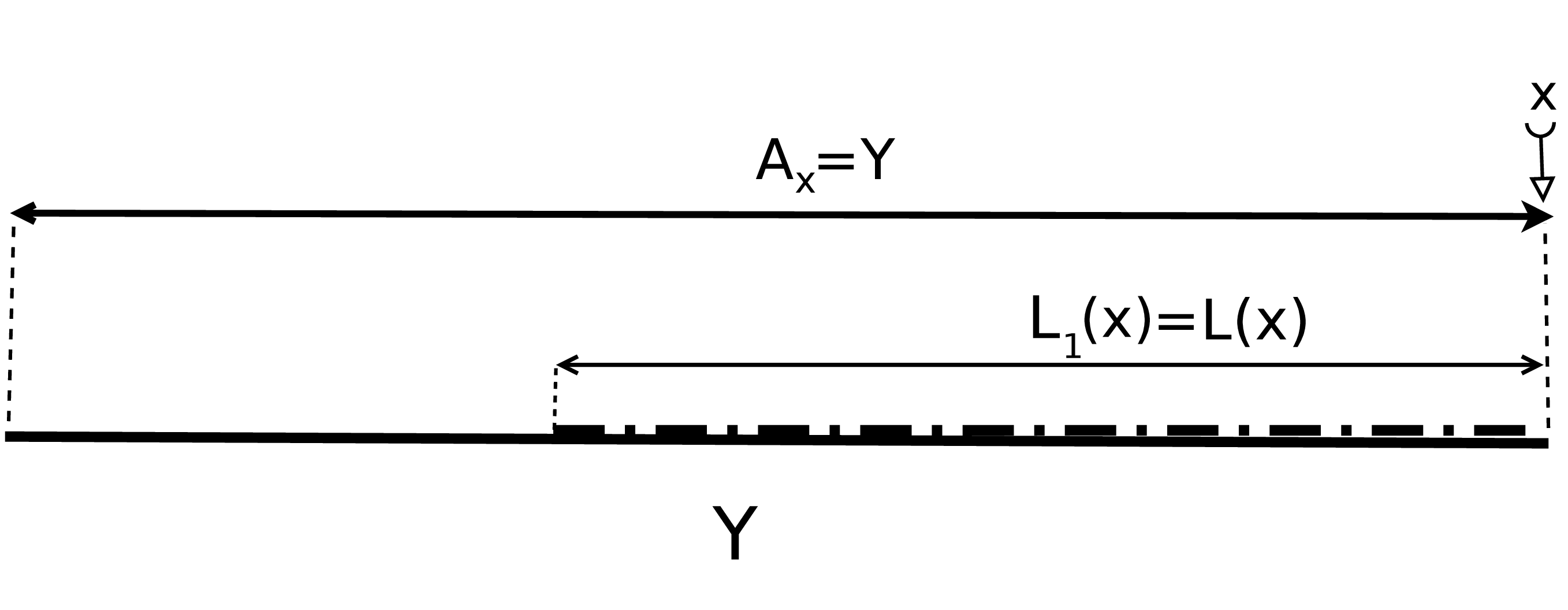}\\
\textbf{Fig.9}.\;\; \begin{small}Points of type $I_c$ such that  $L_1(x)=L(x)$ .                                                                   \end{small}
\end{center}
\medskip
 \item [(iii)] If $L_1(x)=\emptyset$, that is for every 
$y\in L(x)$ and every $y'\in Y$ with $y'<L(x)$ we have that 
$(yxy')_\cam$, then we set
$$\begin{array}{llll}
A_x=L^-(x),\;\;\;B_x=L(x).
\end{array}$$ 

\begin{center}
\includegraphics[scale=0.19]{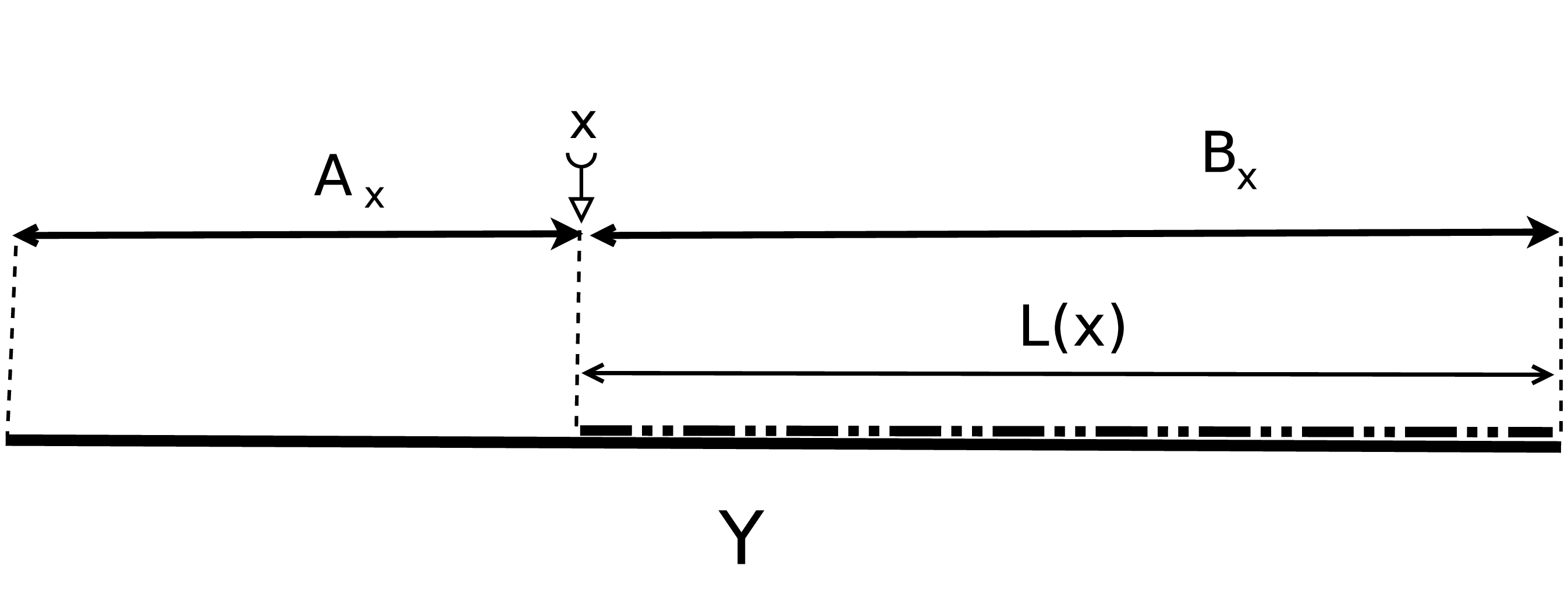}\\
\textbf{Fig.10}.\;\; \begin{small}Points of type $I_c$ such that  $L_1(x)=\emptyset$                                                                        \end{small}.
\end{center}
\medskip
\end{itemize}
\end{definition}
\medskip

\begin{definition}[Sections associated to points of type $\mathrm{I_d}$]\hfill
For every 
$w\in W$ we set $$(A_w, B_w)=(\emptyset, Y).$$
\begin{center}
\includegraphics[scale=0.2]{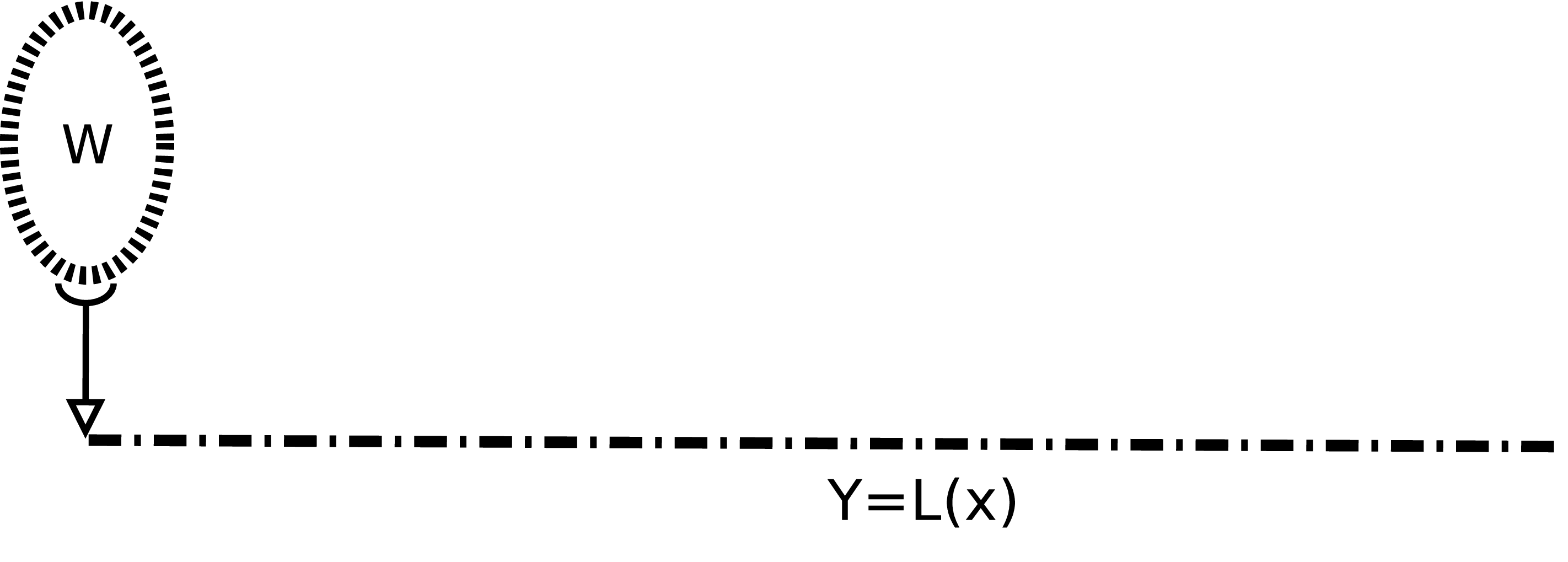}\\
\textbf{Fig.11}\;\;\begin{small} Points of type $I_d$.   \end{small}
\end{center}
\end{definition}

\begin{definition}[Sections associated to points of type $\mathrm{II_a}$]\hfill
\begin{itemize}
  \item[$\mathrm{(v)}$] $ x$ is point of type $\mathrm{II_a}$ and $M(x)=\{y_1, y_2\}$. In that
case we set
$$ A_x=\{ y\in Y: y\leq y_1\}, \;\; B_x=\{y\in Y: y_2\leq y\}.$$
\begin{center}
\includegraphics[scale=0.2]{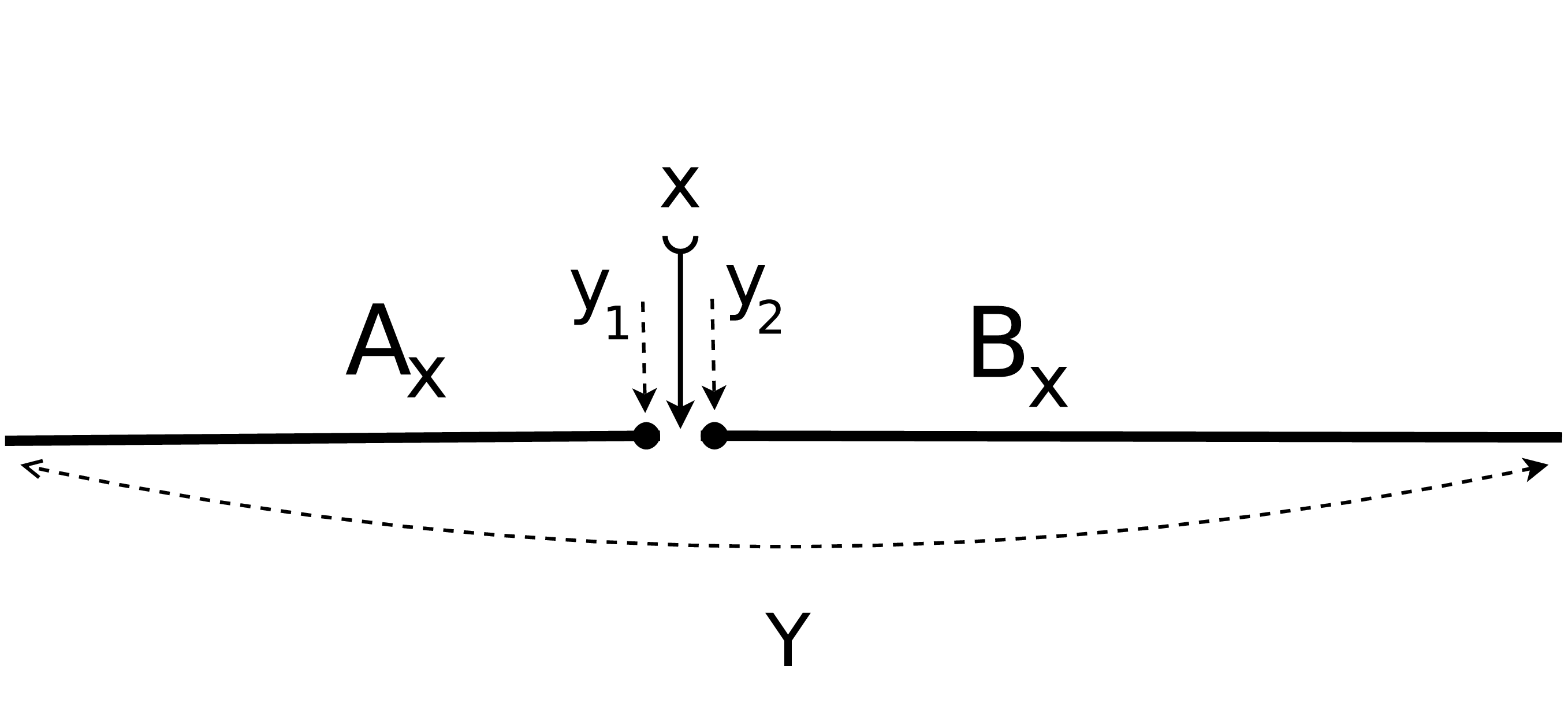}\\
\textbf{Fig.12}\;\; \begin{small}Points $x$ of type $IIa$, with $M(x)=\{y_1, y_2\}$.                                                                  \end{small}
\end{center}
\end{itemize}
\end{definition}

\begin{lemma}\label{lemp11}
If $x$ is a point of type $\mathrm{II_b}$ with $M(x)=\{y_x\}$ then one and only one of the
following cases occurs:
\begin{enumerate}
 \item For every $y>y_x$ we have $\neg(xy_xy)_\cam$.
\item For every $y<y_x$ we have $\neg(xy_xy)_\cam$.
\end{enumerate}
\end{lemma}
 \begin{proof}
Suppose that for some $y>y_x$ we have that $(xy_xy)_\cam$. Let $y'<y_x$
and let $(M', M, A)$ a representative triple for $(y', y, x)$ with $y_x\not \in A$.
Then $A=(M'\vee A)\cap (M\vee A)$  since $y_x\not \in A$ and  $y_x \in M\vee A$
we have that $y_x\not \in M'\vee A$ and so $\neg(xy_xy')_\cam$. The same argument shows that
there exist no $y_1, y_2$  with $y_1<y_x<y_2$ such that $(xy_xy_1)_\cam$ and $(xy_xy_2)_\cam$.
 \end{proof}
\begin{definition}
 If $ x$ is a point of type $\mathrm{II_b}$ and for every $y>y_x$ we have $\neg(xy_xy)_\cam$
we set
$$A_x=\{y\in Y: y\leq y_x\},\;\;  B_x=\{y\in Y: y_x<y\}.$$
\begin{center}
\includegraphics[scale=0.2]{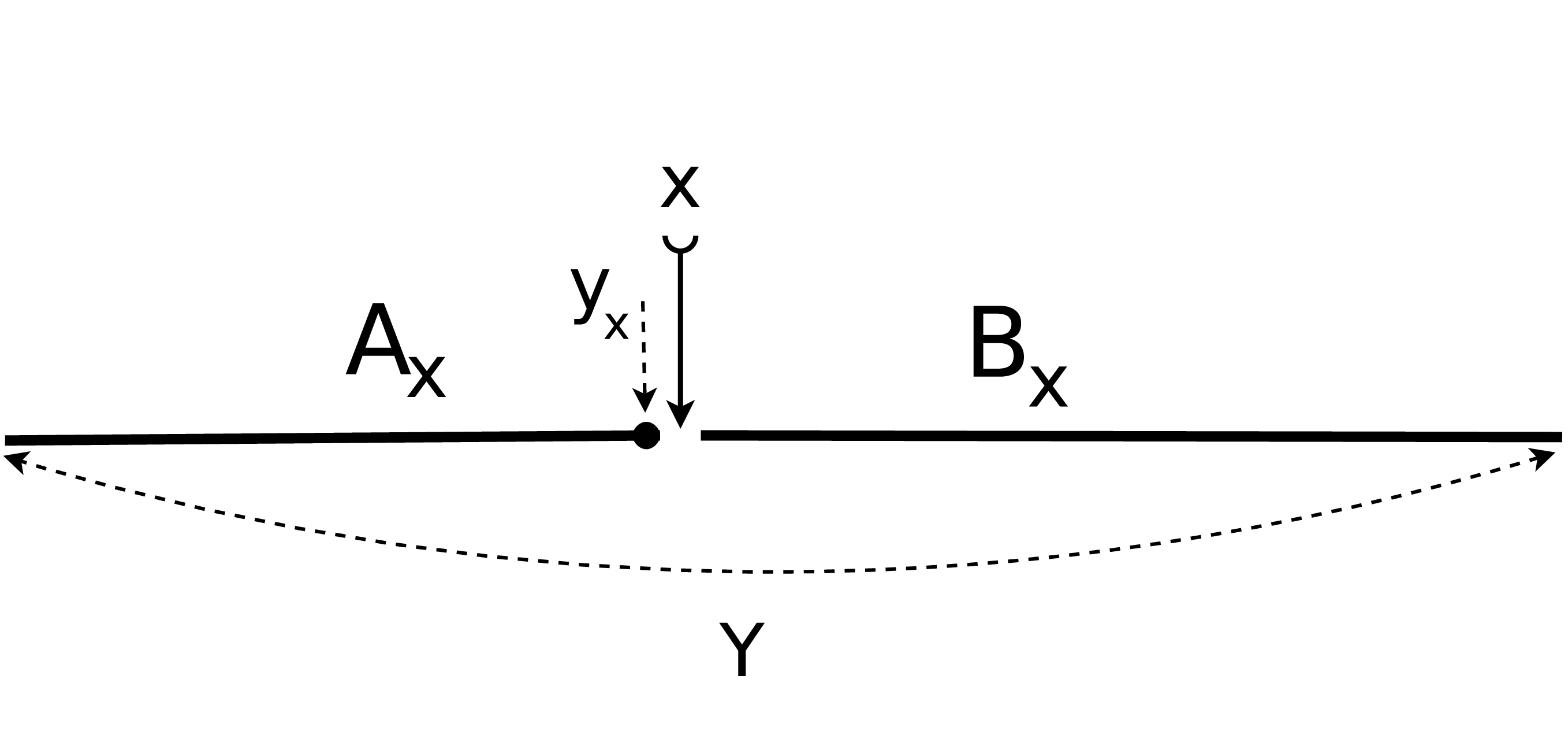}\\
\textbf{Fig.13}\;\; \begin{small}Points $x$ of type $\mathrm{II_b}$ and for every $y>y_x$ we have $\neg(xy_xy)_\cam$.      \end{small}
\end{center}
\end{definition}

\begin{definition}
If  $ x$ is a point of type $\mathrm{II_b}$ and for every $y<y_x$ we have $\neg(xy_xy)_\cam$.
In that case we set
$$A_x=\{y\in Y: y< y_x\},\;\;  B_x=\{y\in Y: y_x\leq y\}.$$
\begin{center}
\includegraphics[scale=0.2]{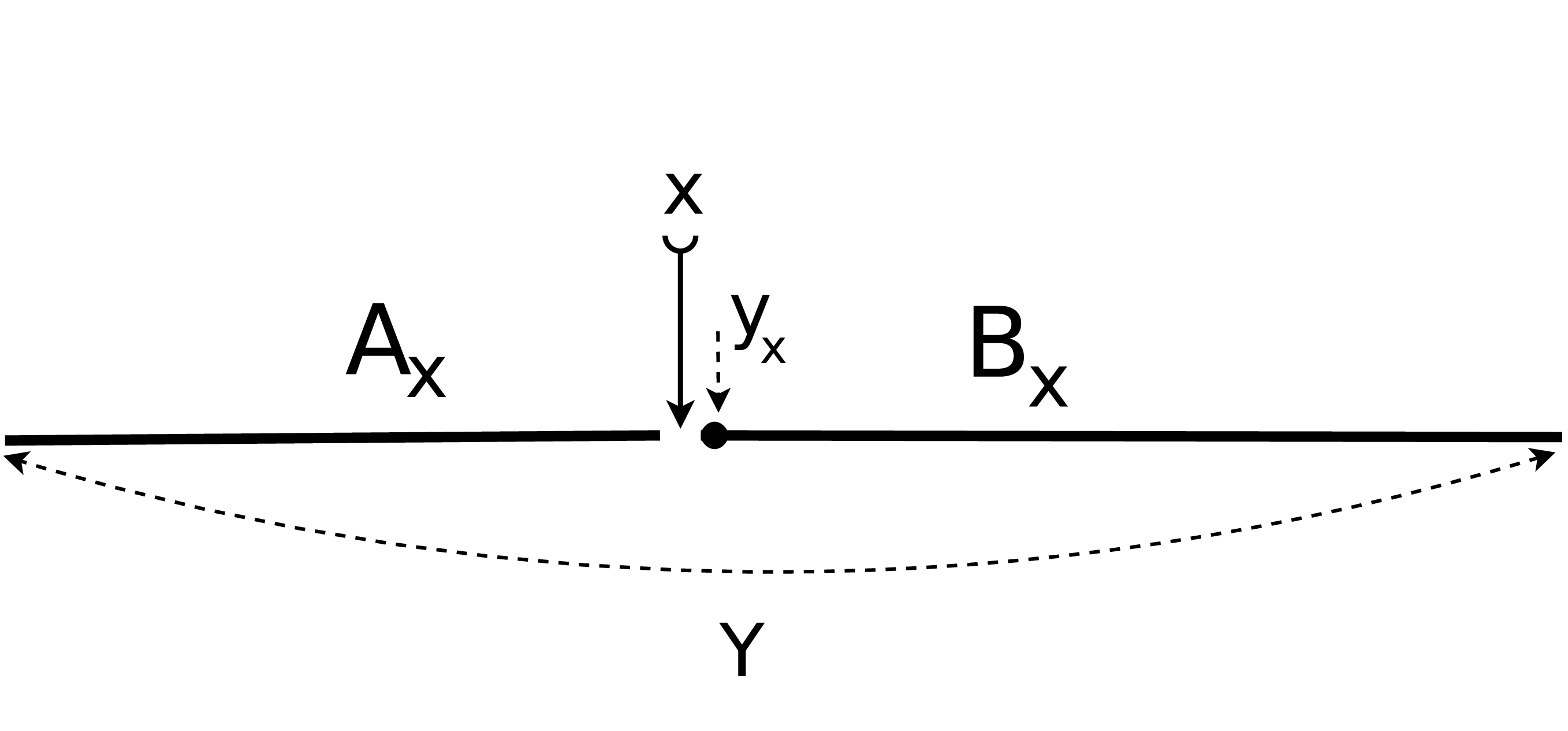}\\
\textbf{Fig.14}\;\; \begin{small}Points $x$  of type $\mathrm{II_b}$ sch that  for every $y<y_x$ we have $\neg(xy_xy)_\cam$.                                                                                           \end{small}   \end{center}

\end{definition}

\begin{definition}
If $ x$ is a point of type $\mathrm{II_c}$, $M(x)=\{y_0\}$
where $y_0$ is  the first element of $Y$  and for some $y>y_0$ we have that $\neg(xy_0y)_\cam$
then we set
$$A_x=\{y_0\},\;\;  B_x=\{y\in Y: y_0<y\}$$\begin{center}
\includegraphics[scale=0.2]{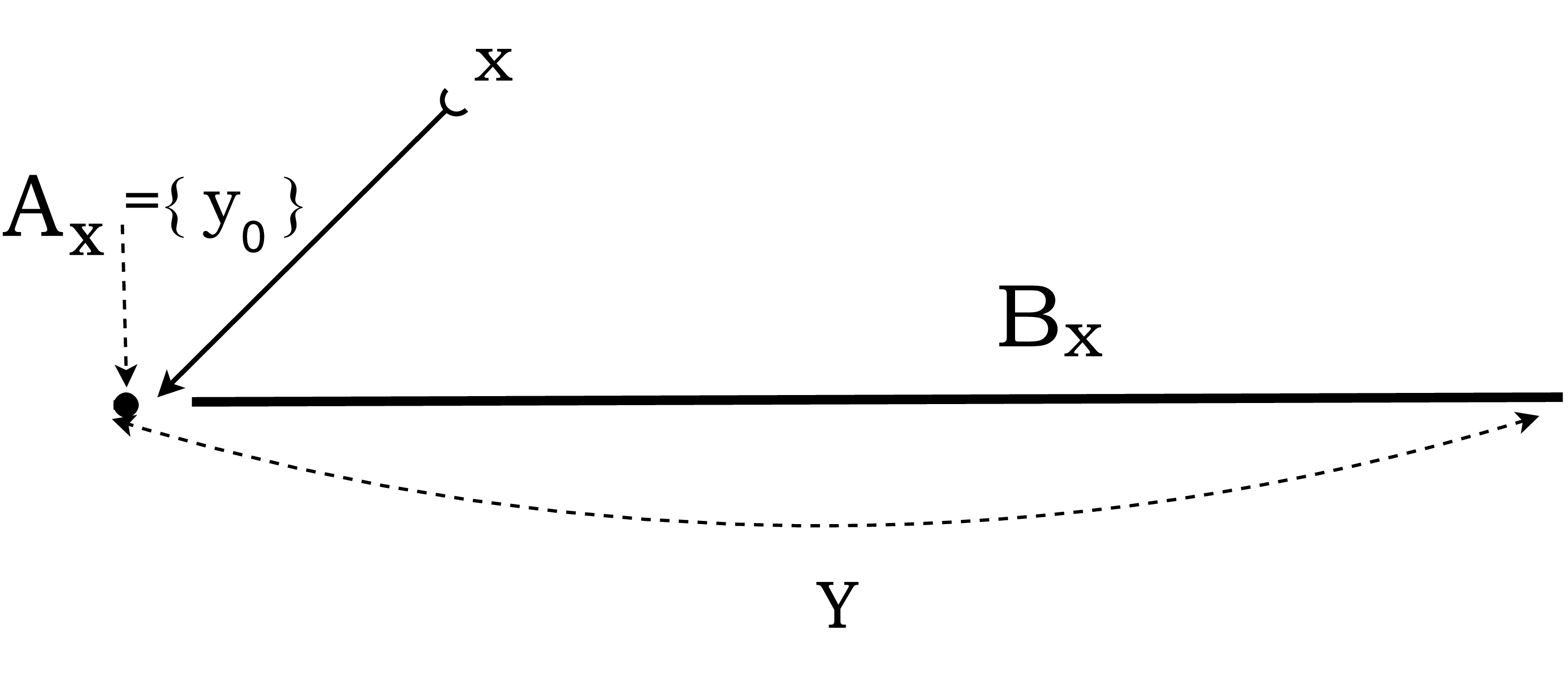}\\
\textbf{Fig.15}\;\; \begin{small}Points $x$  of type $\mathrm{II_c}$
such that  for some $y>y_0$ we have that $\neg(xy_0y)_\cam$.                                        \end{small}   \end{center}\end{definition}

\begin{definition}
If $ x$ is a point of type $\mathrm{II_c}$, $M(x)=\{y_0\}$
where $y_0$ is  the first element of $Y$  and for every  $y>y_0$ we have that $(xy_0y)_\cam$
then we set
$$A_x=\emptyset ,\;\;  B_x=Y. $$
\begin{center}\includegraphics[scale=0.2]{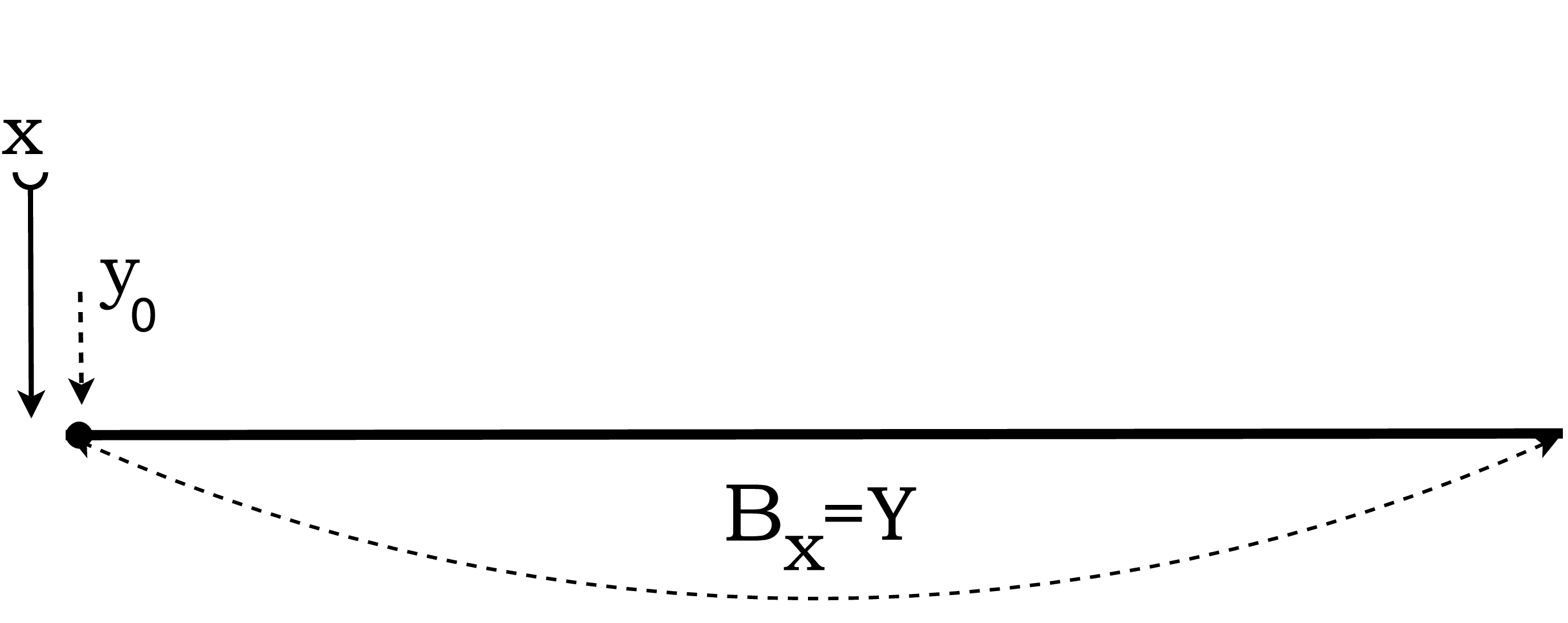}\\
\textbf{Fig.16}\;\; \begin{small}Points $x$  of type $\mathrm{II_c}$
such that  for some $y>y_0$ we have that  $(xy_0y)_\cam$.                                           \end{small}   \end{center}
\end{definition}

\begin{definition}
If $ x$ is a point of type $\mathrm{II_d}$, $M(x)=\{y_1\}$
where $y_1$ is  the last  element of $Y$  and for some $y<y_1$ we have that $\neg(xy_1y)_\cam$
then we set
$$A_x=\{y\in Y: y<y_1\},\;\;  B_x=\{ y_1\}. $$\begin{center}\includegraphics[scale=0.2]{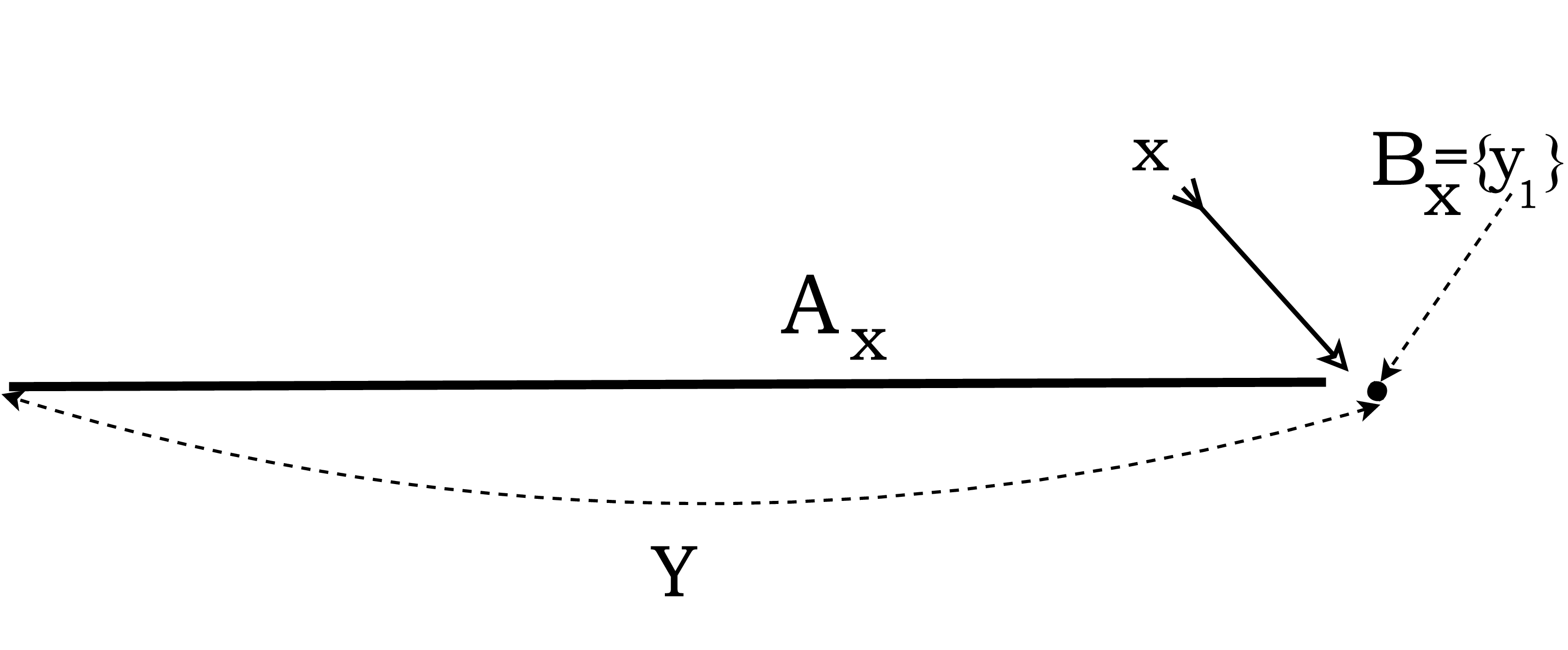}\\
\textbf{Fig.17}\;\; \begin{small}Points $x$  of type $\mathrm{II_c}$
such that  for some $y<y_1$ we have that $\neg(xy_1y)_\cam$.                                           \end{small}   \end{center}\end{definition}
\begin{definition}
If $ x$ is a  point of type $\mathrm{II_d}$, $M(x)=\{y_1\}$
where $y_1$ is  the last  element of $Y$  and for every  $y<y_1$ we have that $(xy_1y)_\cam$
then we set
$$A_x=Y,\;\;  B_x=\emptyset. $$
\begin{center}\includegraphics[scale=0.2]{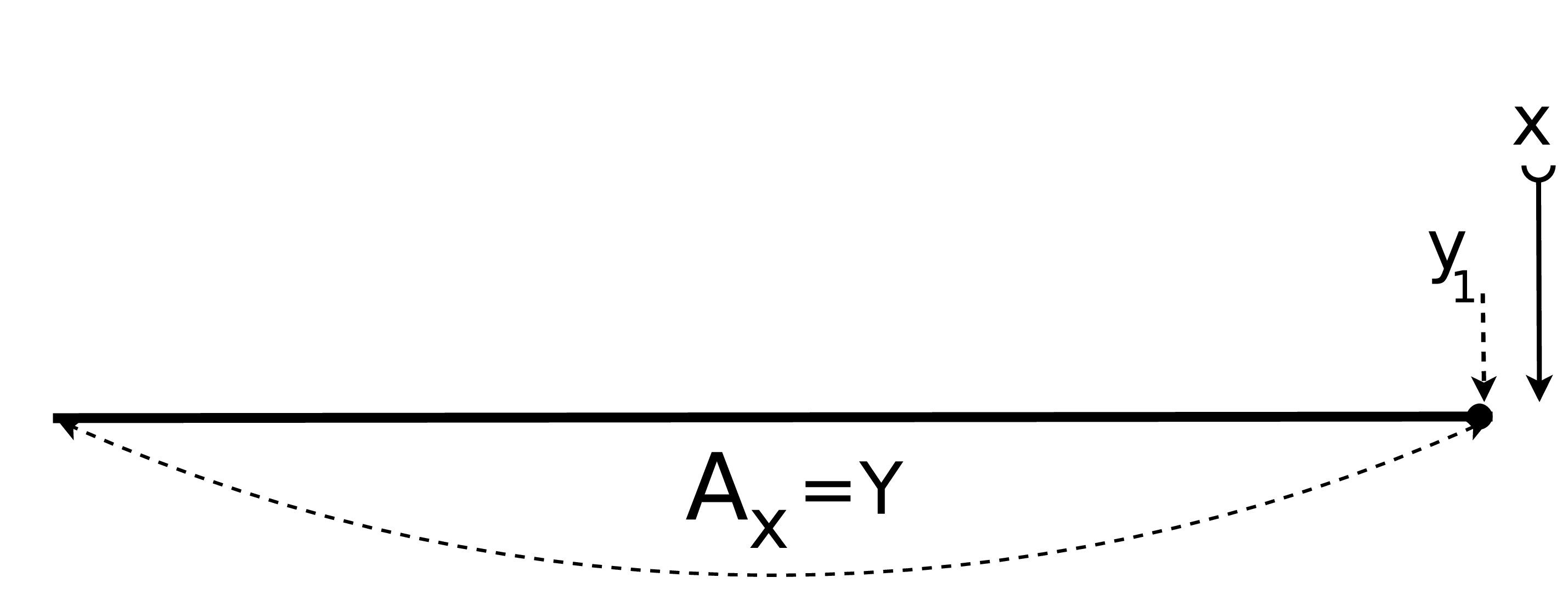}\\
\textbf{Fig.18}\;\; \begin{small}Points $x$  of type $\mathrm{II_c}$
such that  for some $y<y_1$ we have that $(xy_1y)_\cam$.                                           \end{small}   \end{center}\end{definition}

\begin{definition}$ $
\begin{itemize}
 \item [(a)]
We define a partial order $\precsim=\precsim_0$ in $X$ as follows: Let $x, x'\in X$. \\
\begin{itemize}
 \item [(\textbf{i})] If $x, x'\in X\setminus Y$ we set $x\precsim x'$
 if and only if $A_x\subseteq A_{x'}$.\\
\item [(\textbf{ii})]  If $x\in X\setminus Y$ and $x'\in Y$ we set $x\precsim x'$ if and only if $x'\in B_x$.\\
\item [(\textbf{iii})]  If $x\in  Y$ and $x'\in X\setminus Y$ we set $x\precsim x'$ if and only if $x\in A_{x'}$.\\
\item [(\textbf{iv})]  If $x\in  Y$ and $x'\in Y$ we set $x\precsim x'$ if and only if $x\leq x'$.\\
\end{itemize}
\item [(b)]We set $x\sim x'$ if and only if $x\precsim x'$ and $x'\precsim x$.\\
\end{itemize}
\end{definition}

\begin{proposition}
 The quasi order $\precsim$ is $\cam$-consistent.
\end{proposition}

\begin{proof}
	We need to show that the quasi-order $\gtrsim$ is $\mathcal{M}$-consistent, which means that $T_{\gtrsim} \subseteq T_{\mathcal{M}}$. Let $u, v, w \in X$ such that $u \gtrsim v \gtrsim w$. We must prove that $(u v w)_{\mathcal{M}}$ holds. We distinguish four cases based on whether the elements belong to $Y$ or $X \setminus Y$:
	
	\textbf{Case 1:} $u, v, w \in Y$. 
	By Definition 2.40 (iv), the restriction of $\gtrsim$ on $Y$ coincides with the linear ordering $\ge$ of $Y$. Since $\ge$ was chosen to be an $\mathcal{M}$-consistent linear ordering of the maximal independent set $Y$, the relation $(u v w)_{\mathcal{M}}$ holds trivially.
	
	\textbf{Case 2:} Two elements belong to $Y$ and one belongs to $X \setminus Y$. 
	Let $y_1, y_2 \in Y$ and $x \in X \setminus Y$. We assume that $y_1 \gtrsim x \gtrsim y_2$. By parts (ii) and (iii) of Definition 2.40, this is equivalent to $y_1 \in A_x$ and $y_2 \in B_x$. We must show that $(y_1 x y_2)_{\mathcal{M}}$ holds. 
	
	Assume, towards a contradiction, that $\neg(y_1 x y_2)_{\mathcal{M}}$ holds. This means there exists some $M \in \mathcal{M}$ such that $y_1, y_2 \in M$ but $x \notin M$. Note that by the construction of the sections $S_x = (A_x, B_x)$ for all types of $x$, we always have $A_x < B_x$ in $Y$, which implies $y_1 < y_2$. By Lemma 2.9, since $\neg(y_1 x y_2)_{\mathcal{M}}$, we must have either $(x y_1 y_2)_{\mathcal{M}}$ or $(y_1 y_2 x)_{\mathcal{M}}$. We examine the types of $x$:
	
	\textit{Subcase 2.1: $x$ is of Type II.} 
	In this case, $M(x) \neq \emptyset$. By the definition of $A_x$ and $B_x$ for types $II_a, II_b, II_c, II_d$ (Definitions 2.32-2.39), there always exists some $y^* \in M(x)$ such that $y_1 \le y^* \le y_2$. Since $y_1, y_2 \in Y$ and $Y$ is $\mathcal{M}$-consistent, the relation $y_1 \le y^* \le y_2$ implies $(y_1 y^* y_2)_{\mathcal{M}}$. Thus, since $y_1, y_2 \in M$, we must have $y^* \in M$. However, $y^* \in M(x)$ means $y^* L x$, which implies that any set in $\mathcal{M}$ containing $y^*$ must also contain $x$. Therefore, $x \in M$, which is a contradiction.
	
	\textit{Subcase 2.2: $x$ is of Type I.} 
	In this case, $L(x) \neq \emptyset$ and $M(x) = \emptyset$. We analyze the situation where $(x y_1 y_2)_{\mathcal{M}}$ holds (the case for $(y_1 y_2 x)_{\mathcal{M}}$ is completely symmetric). 
	If $(x y_1 y_2)_{\mathcal{M}}$ holds, then for any $N \in \mathcal{M}$, $x, y_2 \in N \implies y_1 \in N$.
	Since $y_2 \in B_x$, for any Type I point, $y_2$ is either in $L(x)$ (specifically $L_2(x)$) or in $L^+(x)$ (specifically $L_2^+(x)$). 
	If $y_2 \in L(x)$, then $x L y_2$. Thus, $x \in N \implies y_2 \in N \implies x, y_2 \in N \implies y_1 \in N$. This means $x L y_1$, so $y_1 \in L(x)$. But if $y_1 \in L(x)$ and $y_2 \in L^+(x)$, the assumption $\neg(y_1 x y_2)_{\mathcal{M}}$ directly places $y_1$ in $L_2(x)$ (by the definition of $L_2(x)$). However, since $y_1 \in A_x$, it must be that $y_1 \in L_1(x)$. This implies $y_1 \in L_1(x) \cap L_2(x)$, which contradicts Lemma 2.27 (2) ($L_1(x) \cap L_2(x) = \emptyset$). 
	The symmetric argument for $(y_1 y_2 x)_{\mathcal{M}}$ forces $y_2 \in L_1(x) \cap L_2(x)$, again yielding a contradiction. Hence, $(y_1 x y_2)_{\mathcal{M}}$ must hold.
The construction of $(Y_0, \leq_0, \precsim_0)$ is completed.\textbf{Case 3:} One element belongs to $Y$ and two elements belong to $X \setminus Y$. 
Let $x_1, x_2 \in X \setminus Y$ and $y \in Y$. We assume that $x_1 \gtrsim y \gtrsim x_2$. By parts (ii) and (iii) of Definition 2.40, this is equivalent to $y \in A_{x_1}$ and $y \in B_{x_2}$. We must show that $(x_1 y x_2)_{\mathcal{M}}$ holds.

Assume, towards a contradiction, that $\neg(x_1 y x_2)_{\mathcal{M}}$ holds. Then there exists $M \in \mathcal{M}$ such that $x_1, x_2 \in M$ and $y \notin M$. 
We first observe that $y$ cannot belong to $L(x_1) \cup M(x_1)$. If $y \in L(x_1)$, then $x_1 L y$, which implies that any set containing $x_1$ contains $y$. Thus $x_1 \in M \implies y \in M$, a contradiction. If $y \in M(x_1)$, then $y L x_1$, so $x_1 \in M \implies y \in M$, again a contradiction. Thus $y \notin L(x_1) \cup M(x_1)$. Since $y \in A_{x_1}$, it follows by the definition of the sections that $y$ lies strictly to the left of the core of $x_1$ (i.e., $y \in L^-(x_1)$ if $x_1$ is of Type I, or $y < M(x_1)$ if $x_1$ is of Type II).

By an identical argument, $y \notin L(x_2) \cup M(x_2)$. Since $y \in B_{x_2}$, it follows that $y$ lies strictly to the right of the core of $x_2$ (i.e., $y \in L^+(x_2)$ or $y > M(x_2)$).
Therefore, $y$ acts as a strict lower bound for the core elements of $x_1$ and a strict upper bound for the core elements of $x_2$ in the linearly ordered set $Y$. This topological separation relative to $Y$ implies that any $M \in \mathcal{M}$ containing both $x_1$ and $x_2$ must contain their intermediate values in $Y$, including $y$, due to the $\mathcal{M}$-consistency of $Y$. Thus $y \in M$, which is a contradiction. Hence $(x_1 y x_2)_{\mathcal{M}}$ holds.

\textbf{Case 4:} All three elements belong to $X \setminus Y$.
Let $x_1, x_2, x_3 \in X \setminus Y$. We assume that $x_1 \gtrsim x_2 \gtrsim x_3$. By Definition 2.40 (i), this is equivalent to the nesting of their left sections: $A_{x_1} \subseteq A_{x_2} \subseteq A_{x_3}$ (and consequently $B_{x_1} \supseteq B_{x_2} \supseteq B_{x_3}$). We must show that $(x_1 x_2 x_3)_{\mathcal{M}}$ holds.

Assume, towards a contradiction, that $\neg(x_1 x_2 x_3)_{\mathcal{M}}$ holds. Then there exists $M \in \mathcal{M}$ such that $x_1, x_3 \in M$ but $x_2 \notin M$.
Since $Y$ is a maximal independent set, the elements $x_1, x_2, x_3$ project onto the sections $(A_{x_i}, B_{x_i})$ in $Y$. The strict inclusion of these sections guarantees that the topological cuts defined by them in $Y$ are ordered. Specifically, the "core" of $x_2$ (i.e., $L(x_2) \cup M(x_2)$) is topologically trapped between the core of $x_1$ and the core of $x_3$. 
By the construction of the quasi-order and the sections (Definitions 2.28-2.39), this nesting ensures that $x_2$ inherits the $\mathcal{M}$-betweenness property from its projection in $Y$. Any set $M \in \mathcal{M}$ containing $x_1$ and $x_3$ covers the intermediate section where the core of $x_2$ resides. Since $x_2 L y'$ or $y' L x_2$ for any $y'$ in its core, $M$ is forced to contain $x_2$, yielding a contradiction.
Hence $(x_1 x_2 x_3)_{\mathcal{M}}$ holds.

In all cases, the relation $(u v w)_{\mathcal{M}}$ is satisfied, which proves that the quasi-order $\gtrsim$ is $\mathcal{M}$-consistent.
\end{proof}
\medskip

\medskip

\noindent \underline{\textbf{Step $\alpha$}}

\medskip

 Let $\alpha$ be an ordinal  with $\alpha>0$
and suppose that  for any ordinal $\xi<\alpha$  
we have defined $(Y_\xi, \leq_\xi, \precsim_\xi)$
such that 
\begin{itemize}
\item  The set + $Y_\xi$ is a subset of $X$
and  $\leq_\xi$ is an $\cam$-compatible linear ordering of $Y_\xi$.
\item The quasi order  $\precsim_\xi$ is an $\cam$-compatible
partial ordering 
of $X$ which is an  extension 
of $\leq_\xi$.
 \item If  $\xi\leq \zeta<\alpha$ then
$Y_\xi\subseteq Y_\zeta$. 
\item For every ordinal $\xi<\alpha$, $Y_\xi\neq X$. 
\end{itemize} We consider two cases:

\medskip

\noindent\underline{\textbf{Case 1.\; The ordinal $\alpha$ is a limit ordinal. }}

We set
\begin{itemize}
\item $$Y_\alpha=\bigcup_{\xi<\alpha}Y_\xi.$$
\item The ordering $\leq_\alpha$ is defined by
$$y\leq_\alpha z \;\;\text{if there exists}\; \xi<\alpha \;\text{with}\; y, z\in Y_\xi\;\;\text{and}\;\; y\leq_\xi z.$$
\item The quasi-ordering $\precsim_\alpha$ is defined by
$$x\precsim_\alpha x' \;\;\text{if for every}\; \;\xi<\alpha \;\;\text{we have that}\; \;\; x\precsim_\xi x'.$$
\end{itemize}

\noindent \underline{\textbf{Case 2. There exists an ordinal $\beta$ such that  $\alpha=\beta+1$. }}

\medskip
Clearly, we suppose that $Y_\beta\neq \emptyset$.

We define an equivalence relation $\sim$ on $X$ by 
$$ x\sim y \;\;\text{if and only if}\;\; x\precsim_\beta x'\;\;\text{and}\;\;  x'\precsim_\beta x.$$
Let $$ \tilde X_\alpha=\{[y]_\sim : y\in X\setminus Y_\beta\} $$ be the set of all equivalence classes 
of all points of $X\setminus Y_\beta$. 
Note that if  $y\in X\setminus Y_\beta$ then $[y]_\sim\cap Y_\beta=\emptyset$. 

\no For every $A\in  \tilde X_\alpha$, $A\neq \emptyset$ we set 
$$ \cam_A=\{A\cap M: M\in \cam\}.$$
Then $\cam_A$ is a non empty  loc-lattice of subsets of $A$
and, by Lemma \ref{prop} and  Zorn's Lemma,   there exists 
 a maximal $\cam_A$-independent 
subset $Y_A$ of $A$ and an $\cam_A$-compatible linear 
order $\leq_A$ of $Y_A$. So,  we can  repeat the construction 
in Step $0$ and we find an  $\cam_A$-consistent
quasi-order $\precsim_A$ of $A$ which extents $\leq_A$.  
We choose for every non-empty $A\in  \tilde X_\alpha$
a such  maximal $\cam_A$-independent subset $Y_A$ of $A$,
 an $\cam_A$-compatible linear 
order $\leq_A$ of $Y_A$ and an  $\cam_A$-consistent
quasi-order $\precsim_A$ of $A$ which extents $\leq_A$. 
We set 
$$ Y_\alpha=Y_\beta\cup\left(\bigcup\nolimits_{A \in  \tilde X_\alpha}
Y_A\right).$$
We define a linear ordering $\leq_\alpha$ of $Y_\alpha$ as 
follows: 
\begin{itemize}
 \item If $y, y'\in Y_\beta$ then we set $y\leq_\alpha y'$ if and only if 
$y\leq_\beta y'$.
\item  If $y\in Y_\beta$ and $y'\in Y_A$ for some $A\in  \tilde X_\alpha$
then we set $y\leq_\alpha y'$ if and only if 
$y\precsim_\beta y'$.
\item  If $y\in Y_A$ and $y'\in Y_A$ for some $A\in  \tilde X_\alpha$
then we set $y\leq_\alpha y'$ if and only if 
$y\leq_A y'$.
\item  If $y\in Y_A$ and $y'\in Y_B$ for some $A, B\in  \tilde X_\alpha$
with $A\neq B$ then we set $y\leq_\alpha y'$ if and only if 
$y\precsim_\beta y'$.
\end{itemize}
We define a quasi  ordering $\precsim_\alpha$ of $X$ as 
follows: 
\begin{itemize}
 \item If $x, x'\in Y_\alpha$ then we set $x\precsim_\alpha x'$ if and only if 
$y\leq_\alpha y'$.
\item  If $x, x'\in A$ for some $A\in  \tilde X_\alpha$
then we set $x\precsim_\alpha x'$ if and only if 
$x \precsim_A x'$.
\item  If $x\in A$ and $x'\in B$ for some $A\in  \tilde X_\alpha$,
 $A\in  \tilde X_\alpha$ and $A\neq B$ 
then we set $x\leq_\alpha x'$ if and only if 
$x\precsim_\beta x''$.
\end{itemize}
\section{Semi prime filters and Proof of Theorem \ref{ath}   }
In this final section we  show that every lattice is isomorphic to a lattice $\cam$ of subsets of a set $X$
which is closed under finite intersections and  separates the points of $X$ (Theorem \ref{threp}).
We use   this representation in order to prove Theorem \ref{ath}.

Let $\cal=(L, \preceq)$ be a lattice.  A \textit{filter} of $\cal$
 is a proper subset $F$ of $L$ such that: If  $a\in F$ and $a\preceq b$ then $b\in F$, and if $a, b\in F$
 then $a\wedge b\in F$.  An \textit{ideal } of  $\cal$
 is a proper subset $I$ of $L$ such that: If  $a\in F$ and $b\preceq a$ then $b\in F$, and if $a, b\in F$
 then $a\vee b\in F$.   A \textit{prime filter} is a filter $F$ such that if $a, b\in L$ and
 $a\vee b\in F$ then $a\in F$ or $b\in F$.

\begin{definition}
Let $\cal=(L, \preceq)$ be a lattice and $a, b\in L$. An element  $x\in L$ is
said to be an  \textbf{internal element} of $a, b$
if $x\preceq a\vee b$ and  $(a\wedge x)\vee (b\wedge x)\neq x$.
By  $\mathbf{I}(a, b)$ we shall denote  the set of all
internal elements of $a, b$ and by  $\bar{\mathbf{I}}(a, b)$ the set $\mathbf{I}(a, b)\cup\{a, b\}$.
\end{definition}

\begin{definition}\label{def32}
A filter $F$ of a lattice $\cal=(L, \preceq)$ is called a \textbf{semi prime filter}
if for every $a, b\in L$ with $a\vee b\in F$ we have that $\bar{\mathbf{I}}(a, b)\cap F\neq \emptyset$.
\end{definition}

\begin{definition}
A lattice $\cal=(L, \preceq)$ is said to be
\begin{enumerate}
 \item [(a)] \textbf{ well separated } if for every semi prime filter $F$ of $\cal$ and every
$a\in L\setminus F$ there exists  $b\in F$ such that $a\not \preceq b$ ;
\item[(b)] \textbf{ completely  separated } if for any $F, G$ semi prime filters of $\cal$ with $F\neq G$
we have that $F\not \subseteq G$ and $G\not \subseteq  F$.
\end{enumerate}
\end{definition}

\begin{remark}
 A completely separated lattice is well separated. Indeed, if $F$ is a semi prime filter and
$a\in L\setminus F$ we select a semi prime filter  $G$ with $a\in G$. Since $F\not \subseteq G$
there exists an element  $b\in F\setminus  G$. It is plain that $ a\not \preceq b$ and $b\not \preceq a$.
The lattice $(\cam, \subseteq)$ with $\cam=\{ \emptyset, \{a\}, \{b\}, \{a, b\},\{a, b, c\}\}$
is an example of a well separated  loc-lattice which is not completely separated. Note also
that if $\cal$ is completely separated loc-lattice then $\cal$ is isomorphic to
a completely separated loc-lattice of sets, which by Lemma  \ref{lem02}
is totally separated. We conclude that a loc-lattice $\cal$ is completely separated if and only if
for any $F, G$ distinct semi prime filters of $\cal$ there exist $a\in F$ and  $b\in G$ with $a\wedge b=0$.
\end{remark}

\begin{lemma}[]\label{lem31}
Given a lattice $\mathcal{L}=(L, \preceq)$,  an ideal $I$ of $\cal$
and a filter $F$ of $\cal$ such that $F\cap I=\emptyset$ there
exists a semi prime filter $\widetilde{F}$ such that  $F\subseteq
\widetilde{F}$ and $\widetilde{F}\cap I=\emptyset$.
\end{lemma}
\begin{proof}Let $\mathbf{F}(\cal)$ the set of all filters of $\cal$ and
 $$\mathcal{F}=\{G\in \mathbf{F}(\mathcal{L}):
F\subseteq G \;\text{and}\;\;  G\cap I=\emptyset \}.$$ Let
$\widetilde{F}$ be a maximal element of $(\mathcal{F}, \subseteq)$.
We shall show that $\widetilde{F}$ is a semi prime filter. For every
filter $F$ and $x\in L$ we set
$$F(x)=\{z\in L: \;\text{there exists} \;\; y\in F\; \text{such that}\;\;y\wedge x \preceq z\}.$$

Clearly, $F(x)$ is a filter, $F\subseteq F(x)$ and $x\in F(x)$. Let
$a, b\in L$  such that $a\vee b\in \widetilde{F}$ and $\{a,
b\}\cap\widetilde{F}=\emptyset$. Since $a, b\not \in \widetilde{F}$
and by the maximality of $\widetilde{F}$ we have that
$$I\cap\widetilde{F}(a)\neq \emptyset \;\;\;\text{and} \;\;\;I\cap \widetilde{F}(b)\neq
\emptyset.$$

 Let $$a'\in I\cap\widetilde{F}(a), \;\;\;\;\;\; b'\in I \cap\widetilde{F}(b)$$

Then there exists $x, y \in \widetilde{F}$ such that $x\wedge a
\preceq a'$ and $y\wedge b \preceq b'$. Since $a', b'\in I$ we have
that $x\wedge a, y\wedge b \in I$.

 Let $c=(a\vee
b)\wedge (x\wedge y)$. Clearly $c\in \widetilde{F}$ and $c\preceq
a\vee b$. But,
$$(c\wedge a)\vee (c\wedge b)=(x\wedge y\wedge a)\vee (x\wedge y\wedge b) \in I.$$ So, $c\neq
(c\wedge a)\vee (c\wedge b)$ and therefore $c\in \mathbf{I}(a, b)$.
\end{proof}
\begin{definition}
Let $\cal=(L, \preceq)$ be a lattice and let  $X_\cal$ be the set of all semi prime filters of $\cal$.
We call the mapping
$f: L\rightarrow \wp(X_\cal)$ defined by
$f(x)=\{ F: x\in F\}$ the \textbf{Stone map} of $\cal$. The topology $\tau$ of $X_\cal$
which has as a basis the set $f(L)$ is called the \textbf{Stone topology} of $X_\cal$.
The space $(X_\cal, \tau)$ is a $T_0$ topological space called the \textbf{Stone space} of $\cal$.
\end{definition}

\begin{theorem}\label{threp}
Let $\cal=(L, \preceq )$ be a lattice and $f:L\rightarrow \wp(X_\cal)$ the Stone map of $\cal$.
 Then $f$  is an isomorphism onto $( f(L), \subseteq )$. Moreover
$f(L), \subseteq)$ is a set lattice closed under finite intersections, $f(L)$ separates $X_\cal$
and for every $a, b\in L$ we have that
$$\begin{array}{llll}
&f(a)\wedge f(b)=f(a\wedge b)=f(a)\cap f(b)&\;\;\;\;\;\;(1)
\\[1ex]& f(a)\vee f(b)= f(a\vee
 b)=\bigcup\left\{f(x): x\in \bar{\mathbf{I}}(a, b)\right\}.&\;\;\;\;\;\;(2)
\end{array}
$$
\end{theorem}
\begin{proof}  Clearly, if $a, b\in L$ and $a\preceq b$ then $f(a)\subseteq
f(b)$. Suppose that $f(a)\subseteq f(b)$ and $a\not \preceq b$.
By Lemma \ref{lem31}  there exists an $F\in X_\cal$ such that
$\{x\in L: a\preceq x\} \subseteq F$ and  $F\cap \{x\in L: x\preceq b\} =\emptyset$,
in particular $a\in F$ and $b\not \in F$.  But then $ F\in f(a)\setminus f(b)$, a contradiction.
So, $a\preceq b$ if and only if $f(a)\subseteq f(b)$ and $f$ is an isomorphism onto $f(L)$.
Let  $F, G\in X_\cal$ with $F\neq G$. Then either $F\not \subseteq G$ or
$G\not \subseteq F$. Suppose that $F\not \subseteq G$.  Clearly, the family  $f(L)$ separates $X_\cal$.   Finally the
equations (1), (2) are obvious from that facts that
 for every semi-prime filter $F$ and every points $a, b$  of $L$
 we have that $a\wedge b \in  F$ if and only
  if  $a\in F$ and  $b\in F$
 and that $a\vee b \in  F$ if and only if there exists  $x\in \bar{\mathbf{I}}(a, b) $
 such that  $x\in F$.
 \end{proof}
\begin{remark}
By duality, it is clear that every lattice $\cal=(L, \preceq)$ can also be represented as  lattice of
of subsets of a set $Y_\cal$,  closed under finite unions. The set $Y_\cal$ will be the set
of all semi prime filters of the dual lattice $L^*=(L, \preceq^*)$, where
$\preceq^*$ is the inverse order of $\preceq$.
The elements of $Y_\cal$ are called the \textbf{semi prime ideals} of $\cal$.
\end{remark}

\begin{lemma}\label{lemfinal}
Let $\cal=(L, \preceq )$ be a lattice and let $f:L\rightarrow \wp(X_\cal)$
be  the Stone map of $\cal$.
Then $(f(L), \subseteq)$ is a loc-lattice of sets. Moreover,  if $\cal$ is well separated then
the family $f(L)$ well  separates the set $X_\cal$.
\end{lemma}
\begin{proof}
By Theorem \ref{threp} the Stone mapping $f$ is an isomorphism
so $(f(L), \subseteq)$ is
lattice of subsets of $X_\cal$ and satisfies  the properties (1), (3) and (4) of Definition\ref{def21}.
We show that $f(L) $ satisfies  property (2) of Definition\ref{def21}.
We note that $f(a)=\emptyset$ if and only if $a=0$.
Suppose that for $A=f(a)$, $B=f(b)$  we have that
$A\cap B\neq \emptyset$. Then $a\wedge b\neq 0$ which implies that  $\mathbf{I}(a, b)=\emptyset$.
Indeed, if  $\mathbf{I}(a, b)\neq \emptyset$ there exists  an internal element $c$ of $a, b$.
Then $a, b, c$ must be incomparable and so by the Property (2) of Definition \ref{def21}
we have that $c=(a\vee c)\wedge (b\vee c)$ which implies that  $a\wedge b\preceq c$.
But then $a\wedge c\neq 0$ and $b\wedge c\neq 0$
and by the property (3)  of Definition \ref{def21} we shall have that
 $c=(a\wedge c)\vee (b\wedge c)$, a contradiction.
So,
$$A\vee B=f(a)\vee f(b) =\bigcup _{c\in \bar{\mathbf{I}}(a, b)}f(c)=f(a)\cup f(b)=A\cup B.$$
Suppose that the lattice $\cal$ is well separated and let $A=f(a)\in f(L)$ and $F\not \in A$.
Let   $b\in L$ such that  $b\not \preceq a$ and $b\in F$. Clearly $F\in B=f(b)$ and $B \not\subseteq A$.
Therefore the family $f(L)$ well separates $X_\cal$.
\end{proof}

It is plain that Theorem \ref{ath} follows from Theorem \ref{thlocset} and Lemma \ref{lemfinal}.

\section{some applications in general topology.}
An  \textit{orderable  topological space } is a topological space $(X, \tau)$
such that there exists a linear ordering $\leq$ of $X$ with the property that the open intervals
of $(X, \leq)$ is a base for $\tau$.
A  \textit{weakly orderable  topological space } is a $T_0$ topological space $(X, \tau)$
such that there exists a linear ordering $\leq$ of $X$ with the property $\tau$ has a base
consisting of convex sets.
A \textit{generalized orderable space} or a \textit{suborderable space} (\cite{ch})  is a
Hausdorff  topological space  $(X, \tau)$
such that there exists a linear ordering $\leq$ of $X$ with the property
that $\tau$ has a basis of convex subsets of $X$. It is known (\cite{lu})  that the class of
generalized ordered spaces coincides with the class of subspaces of
linearly ordered topological spaces.
A linear ordering of $\leq$ of $X$ is said  a \textit{ (Dedekind) complete
linear ordering} if every nonempty subset of $X$ with an upper bound has a least upper
bound (supremum).
A \textit{complete orderable  topological space } is a topological space $(X, \tau)$
such that there exists a complete  linear ordering $\leq$ of $X$ with the property that the open intervals
of $(X, \leq)$ is a base for $\tau$.
An  immediate consequence of Theorems  \ref{ath} and \ref{thlocset}   is the following:
\begin{theorem}
\begin{enumerate}
 \item Every loc-lattice is isomorphic to a basis of a weakly orderable space.
\item A $T_0$  topological space $(X, \tau)$ is weakly orderable if and only if has a basis which is a loc- lattice of
subsets of $X$.
\end{enumerate}
\end{theorem}

The problem of characterization of orderable spaces is considered by many authors.
R. L. Moore (\cite{mo}) and A. Wallace  (\cite{wa})  characterized the orderable continua,
S.  Eilenberg (\cite{eil})  characterized the connected orderable spaces,  H. Herrlich (\cite{he})
the countable and the totally disconnected metric orderable spaces.
A   characterization of general orderable and suborderable  spaces  is given by J. van Dalen and E. Wattel (\cite{dw})  and E. De\`ak (\cite{dea}) (see also \cite{fe}, for another approach).
We can obtain further characterizations of orderable and suborderable spaces
using Theorems \ref{ath} and \ref{thlocset}. Before this we shall investigate which elements
of a lattice can be representated as open intervals.

An element $a$ of a lattice $(L, \preceq)$ is said to be
 \textit{accessible from below}
if there exists a subset $A$ of $L$ such that
$\bigvee A=a$ and  $a\not \in A$.  We say that $a$ is \textit{inaccessible from below}
if it is not accessible from below.  Respectively we say that
$a$ is  \textit{accessible from above}
if there exists a  subset $A$ of $L$ such that
$\bigwedge A=a$ and  $a\not \in A$ and  that $a$ is \textit{inaccessible from above}
if it is not accessible from above .

An \textit{open interval} of linear ordered set $(X, \cal)$ is a subset of $X$ of the form
$(a, b)=\{x\in X: a<x<b)$ or $(a, \rightarrow)=\{x\in X: a<x\}$ or
$(\leftarrow, a)=\{x\in X: x<a\}$ or $X$.
\begin{lemma}\label{lem41}
Let $\cal=(L, \preceq)$ be a totally separated loc-lattice and let $a\in L$ be  an inaccessible from above
element of $L$. If $f: L\rightarrow \wp(X_\cal)$ is the Stone map of $\cal$
and $\leq$ the $f(L)$-consistent linear ordering of $X_\cal$ then $f(a)$ is an open interval
\end{lemma}
\begin{proof} Suppose that $X_\cal \neq f(a)$ and that the set
$f(a)^+=\{\mathbf{x}\in X_\cal: f(a)<x\}$ is not empty. We shall show that $f(a)$ has a first
element $\mathbf{b}$. Suppose that $f(a)$ has no a first element. Let
$\cam=\{ f(a) \vee f(x):  f(x)\subseteq f(a)^+$. Then $\bigcap \cam =f(a)$. Indeed,
if there exists a $\mathbf{x}\in \bigcap \cam\setminus f(a)$ then we select
a $\mathbf{x}'\in f(a)^+$ such that $\mathbf{x}'<\mathbf{x}$ and $x, x'\in L$
such that $\mathbf{x}\in f(x)$, $\mathbf{x}'\in f(x')$ and $f(x)\cap f(x')=\emptyset$.
Then $\mathbf{x}\not \in f(a)\vee f(x')$ which contradicts the assumption that $\mathbf{x}
\in \bigcap \cam$.  This implies that $f(a)=\bigcap \{ f(x): f(a)\subseteq f(x), f(a)\neq f(x)\}$.
Since $f$ is an embedding from $(L, \preceq)$ onto $(f(L), \subseteq )$ we conclude that
$a=\bigwedge \{ x\in L: a\prec  x\}$, which contradicts the assumption that
$a$ is inaccessible from above.

 Similarly, if
the set  $f(a)^-=\{\mathbf{x}\in X_\cal: x<  f(a)\}$ is not empty we show that the set
$f(a)^-$ has a last element $\mathbf{a}$. \end{proof}

\begin{lemma}\label{lem42}
Let $\cal=(L, \preceq)$ be a totally separated loc-lattice  such that
every $a\in L$ is  inaccessible from above.  If $f: L\rightarrow \wp(X_\cal)$ is the Stone map of $\cal$
and $\leq$ the $f(L)$-consistent linear ordering of $X_\cal$ then  $(X_\cal, \tau)$ is a complete
linear ordered topological space.
\end{lemma}
\begin{proof}
By Lemma \ref{lem41} for every $a\in L$ we have that $f(a)$ is an open interval. Since
$\cal=(L, \preceq)$ is  a totally separated loc-lattice, then the family $f(L)$ totally separates
the set $X_\cal$ and therefore $f(L)$ is a basis for the topology $\tau_\leq$ generated
by the open intervals of $(X_\cal, \leq)$. So $\tau=\tau_\leq$ and so $(X_\cal, \tau)$ is an ordered space.
It remains to show that $(X, \leq)$ has no gaps. Suppose that $\mathbf{S}=(S_1, S_2)$ is a gap
of $(X_\cal, \leq)$. It is easy to see that the set
$$F_\mathbf{S}=\{a\in L: f(a)\cap S_1\neq \emptyset \;\text{and}\; f(a)\cap S_2\neq \emptyset\}$$
is a semi prime filter of $\cal$ and so an element of $X_\cal$. We set $\mathbf{x}=F_\mathbf{S}$.
Then either $\mathbf{x}\in S_1$ and so $\mathbf{x}$ is the last element of $S_1$
or $\mathbf{x}\in S_2$ and so $\mathbf{x}$ is the first element of $S_2$. Every case contradicts the assumption
that $\mathbf{S}$ is  a gap.\end{proof}

By Lemma \ref{lem42} and Theorems  \ref{ath} and \ref{thlocset} we easily obtain the following:

\begin{theorem}
\begin{enumerate}
 \item Every totally separated loc-lattice  such that
every $a\in L$ is  inaccessible from above is isomorphic to a basis of a complete  orderable space.
\item A $T_1$  topological space $(X, \tau)$ is  orderable if and only if has a basis $\cam$
such that $(\cam, \subseteq)$  is a loc- lattice of subsets of $X$ such that every $M\in \cam$ is  inaccessible from above.
\end{enumerate}
\end{theorem}
\begin{remark}
A family  $\cam$ of sets is called \textbf{interlocking} (see \cite{dw})  provided that every set $M\in \cam$ which is an intersection of strictly larger members of $\cam$ has a representation as a union of strictly smaller members
of $\cam$. So, we may call a lattice $\cal=(L, \preceq)$ to be an \textbf{interlocking lattice}
if every element of $L$  either is inaccessible from above or it is accessible from above and below.
In such lattices the sets $f(a)$, $a\in L$ are either open intervals or unions of open intervals. Therefore
an interlocking and totally separated lattice has a representation as a basis of an ordered space.
\end{remark}
\section*{Acknowledgments}
The authors would like to acknowledge the use of an AI language model for proofreading, LaTeX formatting assistance, and structural refinements during the preparation of this preprint.

\end{document}